\documentclass[leqno,final]{siamltex}
\usepackage{subfig}
\usepackage{graphicx} 
\usepackage{amsmath,amstext,amssymb,bm}
\usepackage{leftidx}

\usepackage{xcolor} 
\usepackage{soul} 
\usepackage{tikz}
\usetikzlibrary{shapes,arrows}
\usepackage{mathrsfs}
\usepackage{epstopdf}
\usepackage{color}
\usepackage{multirow}
\usepackage{tabularx}
\usepackage[shortlabels]{enumitem}

\numberwithin{equation}{section}
\newtheorem{remark}{Remark}[section]

\allowdisplaybreaks[4]

\newcommand{\E}{\mathcal{E}}

\newcommand{\cK}{\mathcal{K}}

\newcommand{\mP}{\mathbb{P}}

\newcommand{\mE}{\mathbb{E}}

\newcommand{\eps}{\epsilon}

\newcommand{\Ome}{\Omega}
\newcommand{\p}{\partial}
\newcommand{\nab}{\nabla}

\def\D{\mathcal{D}}
\def\E{\mathbb{E}}
\def\P{\mathbb{P}}

\begin{document}
	
	\title{Higher order time discretization method for a class of semilinear stochastic partial differential equations with multiplicative noise}
	\markboth{YUKUN LI AND LIET VO AND GUANQIAN WANG}{THE STOCHASTIC PDEs}		
	
	\author{Yukun Li\thanks{Department of Mathematics, University of Central Florida, Orlando, FL, 32816, U.S.A. (yukun.li@ucf.edu). This author was partially supported by the NSF grant DMS-2110728.}
		\and Liet Vo\thanks{Department of Mathematics, Statistics and Computer Science, The University of Illinois at Chicago, Chicago, IL 60607, U.S.A. (lietvo@uic.edu).}
		\and Guanqian Wang\thanks{Department of Mathematics, University of Central Florida, Orlando, FL, 32816, U.S.A. (guanqian.wang@ucf.edu). This author was partially supported by the NSF grant DMS-2110728.}}

	\maketitle
	
	\begin{abstract} In this paper, we consider a new approach for semi-discretization in time and spatial discretization of a class of semi-linear stochastic partial differential equations (SPDEs) with multiplicative noise. The drift term of the SPDEs is only assumed to satisfy a one-sided Lipschitz condition and the diffusion term is assumed to be globally Lipschitz continuous. Our new strategy for time discretization is based on the Milstein method from stochastic differential equations. We use the energy method for its error analysis and show a strong convergence order of nearly $1$ for the approximate solution. The proof is based on new H\"older continuity estimates of the SPDE solution and the nonlinear term. For the general polynomial-type drift term, there are difficulties in deriving even the stability of the numerical solutions. We propose an interpolation-based finite element method for spatial discretization to overcome the difficulties. Then we obtain $H^1$ stability, higher moment $H^1$ stability, $L^2$ stability, and higher moment $L^2$ stability results using numerical and stochastic techniques. The nearly optimal convergence orders in time and space are hence obtained by coupling all previous results. Numerical experiments are presented to implement the proposed numerical scheme and to validate the theoretical results.
	\end{abstract}

	
	
	
	
	\section{Introduction}\label{sec-1}
	We consider the following initial-boundary value problem for general semi-linear stochastic partial differential equations (SPDEs) with function-type multiplicative noise:
	\begin{alignat}{2}
		du &=\bigl[\Delta u + F(u) \bigr]\,dt + G(u) \, dW(t)  &&\qquad\mbox{a.s. in}\,(0,T)\times D,\label{eq1.1}\\
		u &=  0 && \qquad\mbox{a.s. on } (0,T)\times\p D, \label{eq1.2}\\
		u(0)&= u_0 &&\qquad\mbox{a.s. in}\, D,\label{eq1.3}
	\end{alignat}
	where $D = (0, L)^d\subset \mathbb{R}^d \, (d=1,2)$.
	$F, G$ are two given functions that will be specified later. $\{W(t); t\geq 0\}$ denotes an ${{\mathbb{R}}}$-valued Wiener process.
	
	The corresponding stochastic ordinary differential equations of \eqref{eq1.1} (without the Laplacian term) are studied in \cite{kloeden1991numerical,mao2007stochastic} for the case when both $F$ and $G$ are Lipschitz continuous, and in \cite{higham2002strong} for the case when $G$ satisfies the one-sided Lipschitz condition as stated in \eqref{oneside_Lip}. The strong and weak divergence is considered in \cite{hutzenthaler2010strong} for some $F$ which are not Lipschitz continuous. Besides, the corresponding stochastic partial differential equations of \eqref{eq1.1} when $F$ is Lipschitz and non-Lipschitz continuous and when $G$ is additive and multiplicative are studied in \cite{feng2014finite,feng2017finite,feng2021strong,prohl2014strong,majee2018optimal} based on the variational approach and in \cite{brehier2018analysis,gyongy2005discretization,gyongy2016convergence,jentzen2015strong,kovacs2015backward,kovacs2018discretisation,liu2019strong} based on the semigroup approach. Here the half-order convergence is established in \cite{majee2018optimal} when the drift term is $F(u)=u-u^3$ using the Euler-type scheme. The half-order convergence is established in \cite{feng2021strong} for the drift term in \eqref{eq20180812_1} and diffusion term in assumptions {\bf(A1)}--{\bf(A3)} for a fully discrete scheme.
	
	The primary goal of this paper is to design and analyze a first-order numerical scheme for the time discretization of the problem \eqref{eq1.1}--\eqref{eq1.3}. Specifically, we design a new time discretization method first and then propose an interpolation finite element method, which is based on the new time scheme to discretize the space. Our idea for the time discretization method is inspired by the Milstein method \cite{mil1975approximate} from stochastic differential equations and the semi-discrete in time strategy of the stochastic Stokes equations in \cite{vo2022higher}. In addition, the diffusion function $G$ is assumed to satisfy the global Lipschitz condition while the drift-nonlinear function $F$ is only one-sided Lipschitz. Furthermore, to establish the rates of convergence of the proposed scheme, we use the energy method followed by two steps: the first step is to prove the first-order error order in time by utilizing several established H\"older continuity estimates. The second step is to prove the optimal error order in space. To achieve this, the $H^1$ stability of the numerical solution is needed. The $H^1$-seminorm stability of the numerical solution is proved first and based on which the $L^2$ stability of the numerical solution is established.
	
	The remainder of this paper is organized as follows. In Section \ref{sec2}, several H\"older continuity results about the strong solution are proved. These results will be used in establishing the semi-discrete in-time error estimates. In Section \ref{section_semi}, we present the new approach for the time discretization and its a priori stability as well as the error estimates of the semi-discrete solution are proved. The convergence order is proved to be nearly $1$ for the proposed scheme in $L^2$-norm and the energy norm. In Section \ref{fully_discrete}, we consider an interpolation finite element method for spatial discretization. The finite element method is designed where the interpolation operator is utilized to overcome the difficulty resulting from nonlinearity. Through this approach, the second moment and higher moment $H^1$ stability results are proved first, based on which the second moment and higher moment $L^2$ stability results are proved. Finally, the error estimates with optimal convergence order in space are established based on those stability results. In Section \ref{nume}, several numerical tests including different initial conditions, drift terms, and diffusion terms are used to validate the theoretical results.
	
	\section{Preliminaries}\label{sec2}
	Let $\mathcal{T}_h$
	be the triangulation of $\D$ satisfying the following assumption \cite{xu1999monotone}:
	\begin{equation}\label{eq20180907}
		\frac{1}{d(d-1)}\sum_{K\supset E}|\kappa_E^K|\cot\theta_E^K\ge0,
	\end{equation}
	where $E$ denotes the edge of simplex $K$. 
	Note this assumption is just the Delaunay triangulation when $d=2$. In 3D, the notations in the assumption \eqref{eq20180907} are as follows: $a_i\ (1\leq i \leq d+1)$ denote the
	vertices of $K$, $E=E_{ij}$ the edge connecting two vertices $a_i$ and
	$a_j$, $F_i$ the $(d-1)$-dimensional simplex opposite to the vertex
	$a_i$, $\theta_{ij}^K$ or $\theta_E^K$ the angle between the faces
	$F_i$ and $F_j$, and $\kappa_E^K=F_i \cap F_j$.
	
	Let ${\mathcal{H}}$, ${\cK}$ be two Hilbert spaces. Then, $\mathcal{L}({\mathcal{H}},{\cK})$ is the space of linear maps from ${\mathcal{H}}$ to ${\cK}$. For $m \in \mathbb{N}$, inductively define
	\begin{align}
		\mathcal{L}_m({\mathcal{H}}, {\cK}) := \mathcal{L}({\mathcal{H}}, \mathcal{L}_{m-1}({\mathcal{H}},{\cK})),
	\end{align} 
	as the space of all multi-linear maps from ${\mathcal{H}}\times\cdots\times {\mathcal{H}}$ ($m$ times) to ${\cK}$ for $m \geq 2$.
	
	For some function $G: {\mathcal{H}} \rightarrow {\cK}$, we define the Gateaux derivative of $G$ with respect to $u \in {\mathcal{H}}$, $DG(u) \in \mathcal{L}({\mathcal{H}},{\cK})$, whose action is seen as
	\begin{align*}
		v \mapsto DG(u)(v)\qquad\forall v \in {\mathcal{H}}.
	\end{align*}
	
	In general, we denote $D^k G(u)\in \mathcal{L}_m({\mathcal{H}},{\cK})$, as the $k$-Gateaux derivative of $G$ with respect to $u\in {\mathcal{H}}$.
	
	Below, we state the assumptions on the functionals $G, F: {\mathcal{H}} \rightarrow {\cK}$.
	\medskip
	\begin{enumerate}
		\item[{\bf(A1)}]  $G$ is globally
		Lipschitz continuous and has linear growth. Namely, 
		there exists a constant $C > 0$ such that for all ${v}, {w} \in {\mathcal{H}}$ 
		\begin{subequations}\label{G}
			\begin{align}\label{Lip}
				\|G({ v})- G({ w})\|_{{\cK}} &\leq C\|{v}-{w}\|_{{\mathcal{H}}}\, , \\
				\|G({v})\|_{{\cK}}  &\leq C \bigl( \|{ v}\|_{{\mathcal{H}}}+1\bigr)\, .   \label{lineargrow}
			\end{align}
	\end{subequations}
	\item[{\bf (A2)}] There exists a constant $C>0$ such that
	\begin{align}
		\|DG\|_{L^{\infty}({\mathcal{H}};\mathcal{L}({\mathcal{H}},{\cK}) )} +  \|D^2G\|_{L^{\infty}({\mathcal{H}};\mathcal{L}_2({\mathcal{H}},{\cK}) )}  \leq C.
	\end{align} 
	\item[{\bf (A3)}] There exists a constant $C>0$ such that for all $u,v\in \mathcal{H}$
	\begin{align}
		\|(DG(u) - DG(v)){G}(v)\|_{\mathcal{K}} \leq C\|u - v\|_{\mathcal{H}}.
	\end{align}
\end{enumerate}

In this paper, suppose that $G: H^1_0(D)\rightarrow H^1_0(D)$, and
\begin{align}\label{eq20180812_1}
	F(u)=c_0u-c_1u^3-c_2u^5-c_3u^7-\cdots,
\end{align}
where $c_i\ge 0, i=0,1,2,\cdots$. For simplicity, we choose $F(u) = u - u^q$ for all odd numbers $q \geq 3$. Then $F$ satisfies the following one-sided Lipschitz condition \cite{Higham}
\begin{align}\label{oneside_Lip}
	\langle a - b, F(a) - F(b)\rangle \leq \mu |a - b|^2\qquad\forall a, b \in \mathbb{R}^d,
\end{align}
where $\mu$ is a positive constant. 


Under the above assumptions for the drift term and the diffusion term, it can be proved
in \cite{gess2012strong} that there exists a unique strong variational solution u such that
\begin{align}\label{weak_form}
	\bigl(u(t), \phi\bigr) &= \bigl(u(0),\phi\bigr) - \int_0^t \bigl(\nab u(s), \nab \phi\bigr)\, ds \\\nonumber
	\qquad\qquad&+ \int_0^t \bigl( F(u(s)), \phi\bigr)\, ds + \int_0^t \bigl(G(u(s)), \phi\bigr)\, dW(s) \quad\forall \phi  \in H_0^1(D)
\end{align}
holds $\mathbb{P}$-almost surely. Moreover, when the initial condition $u_0$ is sufficiently smooth, the following stability estimate for the strong solution $u$ holds
\begin{align}\label{pde_estimate}
	\sup_{t \in [0,T]} \mE\bigl[\|u(t)\|^{2}_{H^2}\bigr] + \sup_{t \in [0,T]}\mE\bigl[\|u(t)\|^{4q - 2}_{L^{4q -2}}\bigr] \leq C,
\end{align}
where $q$ is the exponent in the drift term of $F(u) = u  - u^q$.

\bigskip

Next, we introduce the H\"older continuity estimates for the variational solution $u$. 

\begin{lemma}\label{lemma2.4}
	Suppose that the solution $u$ of \eqref{weak_form} satisfies \eqref{pde_estimate}. For $\eps >0$, let $\theta_1 = \frac12 - \epsilon >0$, $\theta_2 = 1-\epsilon>0$. There exists a constant $C \equiv C(D, T, q, u_0)>0$, such that for all $s,t \in [0,T]$,
	\smallskip
	\begin{enumerate}
		\item[(i)] 	 $\displaystyle{\mathbb E}\bigl[\| {u}(t)-{ u}(s)\|^{2}_{H^1} \bigr] 
		\leq C|t-s|^{2\theta_1}.$\\
		\item[(ii)] $\displaystyle{\mathbb E}\Bigl[\Bigl\| {u}(t)-{u}(s) - \int_s^t  G(u(\xi))\, d W(\xi)\Bigr\|^{q}_{H^1} \Bigr] 
		\leq C|t-s|^{q\theta_2},$ where $q = 2, 4.$
		\item [(iii)]$\displaystyle{\mathbb E}\bigl[\| {u}(t)-{ u}(s)\|^{q}_{L^q} \bigr] 
		\leq C|t-s|^{q\theta_1},$
		where $q \geq 2$ are integers.
		\item [(iv)] $\displaystyle{\mathbb E}\Bigl[\Bigl\| F({u}(t))-F({u}(s)) - \int_s^t  DF(u(s))G(u(\xi))\, d W(\xi)\Bigr\|^{2}_{L^2} \Bigr] 
		\leq C|t-s|^{2\theta_2}.$
	\end{enumerate}
\end{lemma}
\medskip
\begin{proof}
	The proof of $(i)$ can be found in \cite[Lemma 2.1]{feng2021strong}, while the establishment of $(iii)$ is based on the semigroup theory, which can be found in many references such as \cite{Printems, Lord, LV2021}. In addition, the proof of $(ii)$ is followed \cite[Lemma 10.27]{Lord} and \cite[Lemma 2.3]{vo2022higher} with minor modifications for $q = 4$. We just need to prove $(iv)$. To prove $(iv)$, we use the Taylor expansion for $F$ with respect to $u(s) \in L^2(D)$  as follows.
	\begin{align}
		F(u(t)) = F(u(s)) + DF(u(s))\bigl(u(t) - u(s)\bigr) + R_2,
	\end{align}
	where $\displaystyle R_2 = \int_0^1 (1-\eta) \bigl(D^2F(u(s) + \eta(u(t) - u(s)))\bigr)(u(t) - u(s))^2\, d\eta$.
	
	Therefore, we have
	\begin{align*}
		&F(u(t)) - F(u(s)) - \int_s^t DF(u(s)) G(u(\xi)) \, dW(\xi) \\\nonumber
		&= DF(u(s)) \bigg[u(t) - u(s) - \int_s^t G(u(\xi))\, dW(\xi)\bigg] + R_2.
	\end{align*}
	
	Since we have $DF(u) = 1 - q u^{q-1}$, then we obtain
	\begin{align}\label{eq2.11}
		&\Bigl\|DF(u(s)) \Bigl[u(t) - u(s) - \int_s^t G(u(\xi))\, dW(\xi)\Bigr]\Bigr\|^2_{L^2} \\\nonumber
		&= \int_D \Bigl| (1 - q u(s)^{q-1})\Bigl[u(t) - u(s) - \int_s^t G(u(\xi))\, dW(\xi)\Bigr]\Bigr|^2\, dx \\\nonumber
		&\leq \int_D  2(1 + q^2 |u(s)|^{2(q-1)})\Bigl|\Bigl[u(t) - u(s) - \int_s^t G(u(\xi))\, dW(\xi)\Bigr]\Bigr|^2\, dx \\\nonumber
		&\leq 2 \Bigl(\int_D(1+ q^2|u(s)|^{2(q-1)})^2\, dx\Bigr)^{\frac12} \Bigl\|u(t) - u(s) - \int_s^t G(u(\xi))\, dW(\xi)\Bigr\|^2_{L^4}\\\nonumber
		&\leq 2 \Bigl(\int_D2(1+ q^4|u(s)|^{4(q-1)})\, dx\Bigr)^{\frac12} \Bigl\|u(t) - u(s) - \int_s^t G(u(\xi))\, dW(\xi)\Bigr\|^2_{L^4}.
	\end{align}
	Taking the expectation $\mE[\cdot]$ to \eqref{eq2.11} and then using the Cauchy-Schwarz inequality, we obtain
	\begin{align}
		&\mE\Bigl[\Bigl\|DF(u(s)) \Bigl[u(t) - u(s) - \int_s^t G(u(\xi))\, dW(\xi)\Bigr]\Bigr\|^2_{L^2}\Bigr] \\\nonumber
		&\leq \mE\Bigl[2 \Bigl(\int_D2(1+ q^4|u(s)|^{4(q-1)})\, dx\Bigr)^{\frac12} \Bigl\|u(t) - u(s) - \int_s^t G(u(\xi))\, dW(\xi)\Bigr\|^2_{L^4}\Bigr]\\\nonumber
		&\leq C_q \Bigl(\mE\Bigl[\|u(s)\|^{4(q-1)}_{L^{4(q-1)}}\Bigr]\Bigr)^{\frac12} \Bigl(\mE\Bigl[\Bigl\|u(t) - u(s) - \int_s^t G(u(\xi))\, dW(\xi) \Bigr\|^4_{L^4}\Bigr]\Bigr)^{\frac12}.
	\end{align}
	Using the interpolation inequality that $\mE[\|u\|^4_{L^4}] \leq C \mE[\|u\|^2_{L^2}\|\nab u\|^2_{L^2}] \leq C \mE[\|u\|^4_{H^1}]$ and Lemma \ref{lemma2.4} $(iii)$ yield to
	\begin{align}
		&\mE\Bigl[\Bigl\|u(t) - u(s) - \int_s^t G(u(\xi))\, dW(\xi) \Bigr\|^4_{L^4}\Bigr]\\\nonumber
		&\leq C\mE\Bigl[\Bigl\|u(t) - u(s) - \int_s^t G(u(\xi))\, dW(\xi) \Bigr\|^4_{H^1}\Bigr] \leq C |t - s| ^{4\theta_2}.
	\end{align}
	By using \eqref{pde_estimate}, we arrive at
	\begin{align}\label{eq2.144}
		&\mE\Bigl[\Bigl\|DF(u(s)) \Bigl[u(t) - u(s) - \int_s^t G(u(\xi))\, dW(\xi)\Bigr]\Bigr\|^2_{L^2}\Bigr] \leq C|t - s|^{2\theta_2},
	\end{align}
	where $C = C_q \Bigl(\sup_{s \in [0,T]}\mE\Bigl[\|u(s)\|^{4(q-1)}_{L^{4(q-1)}}\Bigr]\Bigr)^{\frac12}$.
	
	It is remaining to estimate $R_2$. To do that, we notice that $D^2F(u) = -q(q-1)u^{q-2}$. In the end, we have
	\begin{align}\label{eq2.15}
		&\|R_2\|^2_{L^2} \\ \nonumber
		\leq& \int_D\Bigl| \int_0^1 (1-\eta) q(1-q) (u(s) + \eta (u(t) - u(s)))^{q-2} (u(t) - u(s))^2\, d\eta\Bigr|^2\, dx\\\nonumber
		\leq& \int_D\Bigl(q(q-1)2^{q-2} \bigl(|u(s)|^{q-2} + |u(t) - u(s)|^{q-2}\bigr) \Bigr)^2|u(t) - u(s)|^4\, dx\\\nonumber
		\leq& \int_D q^2 (q-1)^2 2^{2q-3} \bigl(|u(s)|^{2(q-2)} + |u(t) - u(s)|^{2(q-2)}\bigr) |u(t) - u(s)|^4\, dx\\\nonumber
		=& C_q \int_D |u(s)|^{2(q-2)} |u(t) - u(s)|^{4}\, dx + C_q \int_D|u(t) - u(s)|^{2q}\,dx\\\nonumber
		\leq& C_q \|u(s)\|^{2(q-2)}_{L^{4(q-2)}} \|u(t) - u(s)\|^4_{L^8} + C_q \|u(t) - u(s)\|^{2q}_{L^{2q}}.
	\end{align}
	Taking the expectation $\mE[\cdot]$ to \eqref{eq2.15}, using Lemma \ref{lemma2.4} $(iii)$ and then \eqref{pde_estimate} , we obtain
	\begin{align}\label{eq2.16}
		\mE[\|R_2\|^2_{L^2}] &\leq C_q\mE\bigl[\|u(s)\|^{2(q-2)}_{L^{4(q-2)}} \|u(t) - u(s)\|^4_{L^8}\bigr] + C_q \mE\bigl[\|u(t) - u(s)\|^{2q}_{L^{2q}}\bigr] \\\nonumber
		&\leq C_q\Bigl(\mE\bigl[\|u(s)\|^{4(q-2)}_{L^{4(q-2)}}\bigr]\Bigr)^{\frac12} \Bigl(\mE\bigl[\|u(t) - u(s)\|^8_{L^8}\bigr]\Bigr)^{\frac12} \\\nonumber
		&\qquad\qquad+ \mE\bigl[\|u(t) - u(s)\|^{2q}_{L^{2q}}\bigr]\\\nonumber
		&\leq C(|t - s|^{4\theta_1} + |t-s|^{2q\theta_1}) \leq C|t-s|^{4\theta_1},
	\end{align}
	where $C = C_q\Bigl(\sup_{s \in [0,T]}\mE\bigl[\|u(s)\|^{4(q-2)}_{L^{4(q-2)}}\bigr]\Bigr)^{\frac12}$.
	
	The proof is complete by combining \eqref{eq2.144} and \eqref{eq2.16}.
\end{proof}

\section{Semi-discretization in time}\label{section_semi} In this section, we follow the strategy of the Milstein scheme in SDEs to propose a new time discretization method of \eqref{eq1.1}. This approach generates a convergence order of nearly $1$ for the approximate solution. 
\subsection{Formulation of the proposed method}

Let $t_0 < t_1 < \cdots < t_N$ be a uniform mesh of the interval $[0,T]$ 
with the time step size $\tau = \frac{T}{N}$. Note that $t_0 = 0$ and $t_N = T$. 

\smallskip
\noindent
\textbf{Algorithm 1} 

Let $u^0 = u_0$ be a given $H^1_0$-valued random variable. Find $u^{n+1} \in H^1_0(D)$ recursively such that $\mP$-a.s. 
\begin{align}\label{milsteinscheme1}
	\bigl(u^{n+1} - u^n, \phi\bigr) + \tau \bigl(\nab u^{n+1}, \nab\phi\bigr)   &= \tau \bigl(F(u^{n+1}), \phi\bigr) + \bigl(G(u^n)\Delta W_{n} \\\nonumber
	&\qquad+ \frac12 DG(u^n)\,G(u^n)\bigl[(\Delta W_n)^2 - \tau\bigr],\phi\bigr),
\end{align}
for all $\phi \in H^1_0(D)$ and $\Delta W_n = W(t_{n+1}) - W(t_n) \sim \mathcal{N}(0,\tau)$.

{\begin{remark}
		The scheme \eqref{milsteinscheme1} will produce a convergence of order nearly $1$. The difference between \eqref{milsteinscheme1} and the standard Euler-Maruyama method is the discretization of the noise term. While the Euler-type schemes, which establish a convergence order of $\frac12$, contain only the term $G(u^n)\Delta W_n$, the scheme \eqref{milsteinscheme1} adds the extra term $\frac12 DG(u^n)\,G(u^n)\bigl[(\Delta W_n)^2 - \tau]$, which is the key point to obtain a higher convergence order.
\end{remark}}

\bigskip
Next, we define $\displaystyle \mathcal{G}: \mathbb{R}^+\times H_{0}^1(D) \rightarrow L^2(D)$ by
\begin{align}\label{eq3.2}
	\mathcal{G}(s;u) := {G}(u) + DG(u) G(u)\, \int_{t_{n}}^s\,d W(r),\qquad t_n \leq s \leq t_{n+1}.
\end{align}
Then we have
\begin{align*}
	\int_{t_{n}}^{t_{n+1}}\mathcal{G}(s;u^n)\, dW(s) &= G(u^n)\Delta W_{n} + DG(u^n)G(u^n) \int_{t_{n}}^{t_{n+1}}\int_{t_{n}}^s d W(r)\, dW(s)\\\nonumber
	&= G(u^n)\Delta W_{n} + \frac12DG(u^n)G(u^n) \bigl[(\Delta W_{n})^2 - \tau\bigr].
\end{align*}
Therefore, we rewrite \eqref{milsteinscheme1} as follow:
\begin{align}\label{milsteinscheme2}
	\bigl(u^{n+1} - u^n, \phi\bigr) +  \tau \bigl(\nab u^{n+1}, \nab\phi\bigr)   =& \tau\bigl(F(u^{n+1}), \phi\bigr)\\
	&+ 	\int_{t_{n}}^{t_{n+1}}\bigl(\mathcal{G}(s;u^n),\phi\bigr)\, dW(s).\notag
\end{align}

Next, we state the following technical lemma that is used to prove the error estimate results of this paper.
\smallskip
\begin{lemma}\label{lemma3.1} Suppose that $G$ satisfies the assumptions ${\bf(A1), (A2), (A3)}$. Let $u_0\in L^2(\Ome;H^1_0(D)\cap H^2(D))$, there exist constants $C>0$ such that the function $\mathcal{G}$ defined in \eqref{eq3.2} satisfies 
	\smallskip
	\begin{enumerate}[\rm (i)]
		\item  $\displaystyle	\|\mathcal{G}(s;u) - \mathcal{G}(s;v)\|_{L^2} \leq C\|u - v\|_{L^2},\qquad\forall s>0,u,v \in L^2(D)$,
		\smallskip
		\item  $\displaystyle \mE\bigl[\bigl\|G(u(s)) - \mathcal{G}(s;u(t_n))\bigr\|^2_{L^2}\bigr] \leq C|s - t_n|^{2(1-\epsilon)}$, for $t_n\leq s < t_{n+1}$ and $\epsilon>0$.
	\end{enumerate}
\end{lemma}
\smallskip
\begin{proof}
	The Lipschitz continuity of $\mathcal{G}$ in $(i)$ is directly obtained from the assumptions of $G$ while the proof of $(ii)$ can be found in \cite[Lemma 10.36]{Lord} with similar arguments.
	
\end{proof}

Next, we will provide the stability estimates of Algorithm 1 in the following lemma. These stability estimates will be used for the proof of the error estimates of the finite element approximation later.
\smallskip
{
	\begin{lemma}\label{lemma3.3} Let $\{u^n\}$ be the solution of Algorithm 1. Then , there exists a constant $C\equiv C(D, T, u_0, p)$ such that
		
		\begin{enumerate}[(a)]
			\item[(i)] $\displaystyle  \sup_{1 \leq n \leq N} \mE\bigl[\|\nab u^n\|^{2^r}_{L^2}\bigr] + \mE\Bigl[\tau\sum_{n=1}^{N}\|\nab u^n\|^{2^{r} - 2}_{L^2} \|\Delta u^{n}\|^2_{L^2}\Bigr] \leq C$, for any integers $r \geq 1$.
			\item[(ii)] $\displaystyle \sup_{1 \leq n \leq N} \mE\bigl[\|\nab u^n\|^p_{L^2}\bigr] \leq C$, for any integers $p \geq 2$.
		\end{enumerate}
	\end{lemma}
	
}

\begin{proof} We just provide the proof of (i) when $r =1$. When $r \geq 2$, the proof is similar to \cite[Lemma 3.1]{Prohl} with minor modifications. So, we skip it to save space.
	
	To begin, we rewrite \eqref{milsteinscheme1} in the strong form as follow:
	\begin{align}\label{eq2.6}
		u^{n+1} - u^n -\tau\Delta u^{n+1} &= \tau F(u^{n+1}) + G(u^n)\Delta W_n \\\nonumber
		&\qquad+ \frac12 DG(u^n)G(u^n)[(\Delta W_n)^2 - \tau].
	\end{align}
	Testing the equation \eqref{eq2.6} by $-\Delta u^{n+1}$ and then using integration by parts we obtain
	\begin{align}\label{eq3.6}
		&\bigl(\nabla(u^{n+1} -u^n), \nab u^{n+1}\bigr) + \tau \|\Delta u^{n+1}\|^2_{L^2} \\
		&\qquad= -\tau\bigl(F(u^{n+1}), \Delta u^{n+1}\bigr)- \bigl(G(u^n), \Delta u^{n+1}\bigr)\Delta W_{n}\nonumber\\
		&\qquad\quad- \frac12 \bigl(DG(u^n)G(u^n),\Delta u^{n+1}\bigr)\bigl[(\Delta W_n)^2 - \tau\bigr]\nonumber\\
		&\qquad:={\tt I + II + III}.\nonumber
	\end{align}  
	By using the integration by parts, we obtain
	\begin{align} \label{eq_3.7}
		{\tt I} &= -\tau\bigl(u^{n+1}, \Delta u^{n+1}\bigr) + \tau \bigl((u^{n+1})^q,\Delta u^{n+1}\bigr)\\\nonumber
		&=  \tau\|\nab u^{n+1}\|^2_{L^2} - \tau q\bigl((u^{n+1})^{q-1}\nab u^{n+1},\nabla u^{n+1}\bigr)\\\nonumber
		&= \tau\|\nab u^{n+1}\|^2_{L^2} - \tau q \int_D (u^{n+1})^{q-1}|\nab u^{n+1}|^2\, dx \leq \tau \|\nab u^{n+1}\|^2_{L^2}, 
	\end{align}
	where the last inequality of \eqref{eq_3.7} is obtained by using the fact that, for all odd $q\geq 3$,
	$\int_D (u^{n+1})^{q-1}|\nab u^{n+1}|^2\, dx \geq 0$.
	
	To bound {\tt II}, we take the expectation and then use the fact that  $\mE[\Delta W_n] = 0$. Namely,
	\begin{align}
		\mE[{\tt II}] &= -\mE\bigl[\bigl(G(u^n), \Delta(u^{n+1} - u^n)\bigr)\Delta W_n\bigr] - \mE\bigl[\bigl(G(u^n), \Delta u^n\bigr)\Delta W_n\bigr]\\\nonumber
		&= \mE\bigl[\bigl(\nab G(u^n), \nab(u^{n+1} - u^n)\bigr)\Delta W_n\bigr]\\\nonumber
		&\leq C\mE[\|\nab u^{n}\|^2_{L^2}|\Delta W_n|^2] + \frac14\mE\bigl[\|\nab(u^{n+1} - u^n)\|^2_{L^2}\bigr]\\\nonumber
		&= C\tau\mE[\|\nab u^{n}\|^2_{L^2}|] + \frac14\mE\bigl[\|\nab(u^{n+1} - u^n)\|^2_{L^2}\bigr].
	\end{align}
	In addition, by using the Cauchy-Schwarz and the assumptions ${\bf (A1), (A2)}$, we have
	\begin{align}\label{eq2.9}
		\mE[{\tt III}] &\leq \frac{C}{\tau}\mE\bigl[\|DG(u^n)G(u^n)\|^2_{L^2}|(\Delta W_n)^2 - \tau|^2\bigr] + \frac{\tau}{4}\mE\bigl[\|\Delta u^{n+1}\|^2_{L^2}\bigr]\\\nonumber
		&\leq \frac{C}{\tau}\mE\bigl[\|G(u^n)\|^2_{L^2}|(\Delta W_n)^2 - \tau|^2\bigr] + \frac{\tau}{4}\mE\bigl[\|\Delta u^{n+1}\|^2_{L^2}\bigr]\\\nonumber
		&\leq C\tau\mE\bigl[\|\nab u^{n}\|^2_{L^2}\bigr] + \frac{\tau}{4}\mE\bigl[\|\Delta u^{n+1}\|^2_{L^2}\bigr],
	\end{align}
	where the last inequality of \eqref{eq2.9} is obtained by  using the fact that $\mE[|(\Delta W_n)^2 -\tau |^2] \leq C\tau^2$.
	
	Substituting all the estimates from ${\tt I, II, III}$ into \eqref{eq2.6} and absorbing the like-terms from the right side to the left side, we obtain
	\begin{align}\label{eq2.10}
		\frac12\mE\bigl[\|\nab u^{n+1}\|^2_{L^2} &- \|\nab u^n\|^2_{L^2}\bigr] + \frac14\mE\bigl[\|\nab(u^{n+1} - u^n)\|^2_{L^2}\bigr] + \frac{\tau}{2}\mE\bigl[\|\Delta u^{n+1}\|^2_{L^2}\bigr]\\\nonumber
		& \leq C\tau\mE\bigl[\|u^{n+1} - u^n\|^2_{L^2}\bigr] + C\tau\mE\bigl[\|\nab u^{n}\|^2_{L^2}\bigr].
	\end{align}
	Next, applying the summation $\sum_{n=0}^{\ell}$, for any $0\leq \ell <N$, we obtain
	\begin{align}
		&\mE\bigl[\|\nab u^{\ell +1}\|^2_{L^2}\bigr] + \sum_{n=0}^{\ell}\mE\bigl[\|\nab(u^{n+1}-u^n)\|^2_{L^2}\bigr] + \tau\sum_{n=0}^{\ell} \mE\bigl[\|\Delta u^{n+1}\|^2_{L^2}\bigr]\\\nonumber
		&\leq C\tau\sum_{n=0}^{\ell}\mE\bigl[\|\nab u^n\|^2_{L^2}\bigr] + \mE\bigl[\|\nab u_0\|^2_{L^2}\bigr] + C\tau\sum_{n=0}^{\ell} \mE\bigl[\|u^{n+1} - u^n\|^2_{L^2}\bigr].
	\end{align}
	The proof is completed by using Gronwall's inequality. 
	
	Finally, the proof of (ii) is followed by using the result from (i) and H\"older inequality.

\end{proof}

\subsection{Error estimates for Algorithm 1}\label{sec3.2}

In this part, we state the first main result of this paper which establishes an $O(\tau^{1-\eps})$ convergence order for the proposed method.

\begin{theorem}\label{theorem_semi} 
	Let $u$ be the variational solution to \eqref{eq1.1} and $\{u^{n}\}$ be generated by  Algorithm 1. Assume that $G$ satisifies ${\bf (A1), (A2), (A3)}$ and $u_0 \in L^{2}(\Ome; H^1_0(D)\cap H^2(D))$. Suppose that $0<\epsilon<1$, then there exists a constant $C = C(D, T, u_0)>0$ such that 
	\begin{align}\label{eq310}
		\sup_{1\leq n \leq N}\mE\Bigl[\|u(t_n) - u^n\|^2_{L^2}\Bigr] + \mE\bigg[\tau\sum_{n=1}^N\|\nab(u(t_n) - u^n)\|^2_{L^2}\bigg] \leq C\, \tau^{2(1-\epsilon)}.
	\end{align}
\end{theorem}
\begin{proof} Denote $e^n := u(t_n) - u ^n$. Subtracting \eqref{milsteinscheme2} from \eqref{weak_form}, we obtain the following error equation
	\begin{align}\label{eq2.13}
		\bigl(e^{n+1} - e^n, \phi\bigr) + \tau\bigl(\nab e^{n+1}, \nab \phi\bigr) &= \int_{t_{n}}^{t_{n+1}} \bigl(\nab(u(t_{n+1}) - u(s)), \nab \phi\bigr)\, ds \\\nonumber
		&\qquad- \int_{t_{n}}^{t_{n+1}} \bigl(F(u(t_{n+1})) - F(u(s)), \phi\bigr)\, ds\\\nonumber
		&\qquad+  \int_{t_{n}}^{t_{n+1}} \bigl(F(u(t_{n+1})) - F(u^{n+1}), \phi\bigr)\, ds\\\nonumber
		&\qquad+ \int_{t_{n}}^{t_{n+1}} \bigl(G(u(s)) - \mathcal{G}(s;u^n),\phi\bigr)\, dW(s).
	\end{align}
	Now, choosing $\phi = e^{n+1}$ and using the identity $2a(a-b) = a^2 -b^2 + (a-b)^2$, we have
	\begin{align}\label{eq2.14}
		&\frac12\bigl[\|e^{n+1}\|^2_{L^2} - \|e^{n}\|^2_{L^2}\bigr] + \frac12\|e^{n+1} - e^n\|^2_{L^2} + \tau\|\nab e^{n+1}\|^2_{L^2} \\\nonumber
		&= \int_{t_{n}}^{t_{n+1}} \bigl(\nab(u(t_{n+1}) - u(s)), \nab e^{n+1}\bigr)\, ds \\\nonumber
		&\qquad- \int_{t_{n}}^{t_{n+1}} \bigl(F(u(t_{n+1})) - F(u(s)), e^{n+1}\bigr)\, ds\\\nonumber
		&\qquad+  \int_{t_{n}}^{t_{n+1}} \bigl(F(u(t_{n+1})) - F(u^{n+1}), e^{n+1}\bigr)\, ds\\\nonumber
		&\qquad+ \int_{t_{n}}^{t_{n+1}} \bigl(G(u(s)) - \mathcal{G}(s;u^n),e^{n+1}\bigr)\, dW(s)\\\nonumber
		&:= {\tt I + II + III + IV}.
	\end{align}
	Next, we bound the right side of \eqref{eq2.14} as follows.
	
	In order to estimate {\tt I}, we add and subtract $\displaystyle \int_s^{t_{n+1}}\nab G(u(\xi))\, dW(\xi)$ for any $t_n\leq s < t_{n+1}$, as follow.
	\begin{align}
		{\tt I} &= \int_{t_{n}}^{t_{n+1}} \Bigl(\nab\Bigl(u(t_{n+1}) - u(s) - \int_s^{t_{n+1}}G(u(\xi))\, dW(\xi)\Bigr), \nab e^{n+1}\Bigr)\, ds \\\nonumber
		&\qquad+\int_{t_{n}}^{t_{n+1}} \Bigl(\int_s^{t_{n+1}}\nab G(u(\xi))\, dW(\xi),\nab e^{n+1}\Bigr)\, ds\\\nonumber
		&:= {\tt I_1 + I_2}.
	\end{align}
	By using Lemma \ref{lemma2.4} $(ii)$, we obtain
	\begin{align}
		\mE[{\tt I_1}] &\leq \int_{t_{n}}^{t_{n+1}} \mE\Bigl[\Bigl\|u(t_{n+1}) - u(s) - \int_s^{t_{n+1}}G(u(\xi))\, dW(\xi)\Bigr\|^2_{H^1}\Bigr]\, ds\\\nonumber
		&\qquad+ \frac{\tau}{4}\mE\bigl[\|\nab e^{n+1}\|^2_{L^2}\bigr]\\\nonumber
		&\leq C\tau^{1+ 2(1-\eps)} + \frac{\tau}{4}\mE\bigl[\|\nab e^{n+1}\|^2_{L^2}\bigr].
	\end{align}
	Next, by the integration by parts we have
	\begin{align}
		{\tt I_2} &= \int_{t_{n}}^{t_{n+1}} \Bigl(\int_s^{t_{n+1}}\nab G(u(\xi))\, dW(\xi),\nab (e^{n+1} - e^n)\Bigr)\, ds \\\nonumber
		&\qquad+\int_{t_{n}}^{t_{n+1}} \Bigl(\int_s^{t_{n+1}}\nab G(u(\xi))\, dW(\xi),\nab e^{n}\Bigr)\, ds\\\nonumber
		&=-\int_{t_{n}}^{t_{n+1}} \Bigl(\int_s^{t_{n+1}}\Delta G(u(\xi))\, dW(\xi),e^{n+1} - e^n\Bigr)\, ds \\\nonumber
		&\qquad+\int_{t_{n}}^{t_{n+1}} \Bigl(\int_s^{t_{n+1}}\nab G(u(\xi))\, dW(\xi),\nab e^{n}\Bigr)\, ds\\\nonumber
		&:= {\tt I_{2a} + I_{2b}}.
	\end{align}
	
	We note that $\mE[{\tt I_{2b}}] = 0$ due to the martingale property of the It\^o integral. So, it is left to estimate ${\tt I_{2a}}$. 
	By using the H\"older inequality, we obtain
	\begin{align}
		{\tt I_{2a}} &= -\int_{t_n}^{t_{n+1}}\bigg(\int_s^{t_{n+1}}\Delta G(u(\xi))\, dW(\xi), e^{n+1} - e^n\bigg)\, ds\\\nonumber
		&\leq 2\bigg\|\int_{t_n}^{t_{n+1}}\int_s^{t_{n+1}}\Delta G(u(\xi))\,dW(\xi)\,ds\bigg\|^2_{L^2} + \frac18\|e^{n+1} - e^n\|^2_{L^2}\\\nonumber
		& = 2\int_D\bigg|\int_{t_n}^{t_{n+1}}\int_s^{t_{n+1}}\Delta G(u(\xi))\,dW(\xi)\,ds\bigg|^2\, d{\bf x}+ \frac18\|e^{n+1} - e^n\|^2_{L^2}\\\nonumber
		&\leq 2\int_D\bigg(\int_{t_n}^{t_{n+1}}\bigg|\int_s^{t_{n+1}}\Delta G(u(\xi))\,dW(\xi)\bigg|\,ds\bigg)^2\, d{ x}+ \frac18\|e^{n+1} - e^n\|^2_{L^2}\\\nonumber
		&\leq 2\tau\int_D\int_{t_n}^{t_{n+1}}\bigg|\int_s^{t_{n+1}}\Delta G(u(\xi))\, dW(\xi)\bigg|^2\, ds\, d{x}  + \frac18\|e^{n+1} - e^n\|^2_{L^2}\\\nonumber
		&= 2\tau \int_{t_n}^{t_{n+1}}\bigg\|\int_s^{t_{n+1}}\Delta G(u(\xi))\, dW(\xi)\bigg\|^2_{L^2}\, ds + \frac18\|e^{n+1} - e^n\|^2_{L^2}.
	\end{align}
	
	By using the It\^o isometry we have
	\begin{align}
		\mE[{\tt I_2}] = \mE[{\tt I_{2a}}] &\leq C\tau^{3}\sup_{\xi\in [0,T]}\mE[\|u(\xi)\|^2_{H^2}] + \frac18\mE[\|e^{n+1} - e^n\|^2_{L^2}]\\\nonumber
		&\leq C\tau^{3} + \frac18\mE[\|e^{n+1} - e^n\|^2_{L^2}].
	\end{align}
	
	Similarly, we can estimate {\tt II} as follows.
	\begin{align}
		{\tt II} &=- \int_{t_{n}}^{t_{n+1}} \Bigl(F(u(t_{n+1})) - F(u(s)) -\int_s^{t_{n+1}} DF(u(s))G(u(\xi))\, dW(\xi), \\\nonumber
		&\qquad e^{n+1}\Bigr)\, ds-\int_{t_{n}}^{t_{n+1}}\Bigl(\int_s^{t_{n+1}}DF(u(s))G(u(\xi))\, dW(\xi), e^{n+1}\Bigr)\, ds\\\nonumber
		&:= {\tt II_1 + II_2}.
	\end{align}
	
	By using Lemma \ref{lemma2.4} $(iv)$ and Poincar\'e's inequality, we obtain
	\begin{align}
		\mE[{\tt II_1}] &\leq C\tau^{1+2(1-\eps)} + \frac{\tau}{4}\mE\bigl[\|\nab e^{n+1}\|^2_{L^2}\bigr].
	\end{align}
	
	To estimate ${\tt II_2}$, we use the same techniques from estimating ${\tt I_{2}}$ and also use \eqref{pde_estimate}, we obtain
	\begin{align*}
		\mE[{\tt II_2}] &= -\mE\Bigl[\int_{t_{n}}^{t_{n+1}}\Bigl(\int_s^{t_{n+1}}DF(u(s))G(u(\xi))\, dW(\xi), e^{n+1}-e^n\Bigr)\, ds\Bigr]\\\nonumber
		&\leq C\mE\Bigl[\Bigl\|\int_{t_{n}}^{t_{n+1}}\int_{s}^{t_{n+1}}DF(u(s))G(u(\xi))\, dW(\xi)\, ds\Bigr\|^2_{L^2}\Bigr] + \frac{1}{8}\mE\bigl[\|e^{n+1} - e^n\|^2_{L^2}\bigr]\\\nonumber
		&\leq C\mE\Bigl[\int_D\Bigl(\int_{t_{n}}^{t_{n+1}} \Bigl|\int_{s}^{t_{n+1}} DF(u(s))G(u(\xi))\, dW(\xi)\Bigr|\, ds\Bigr)^2\, dx\Bigr]\\\nonumber
		&\qquad+ \frac{1}{8}\mE\bigl[\|e^{n+1} - e^n\|^2_{L^2}\bigr]\\\nonumber
		&\leq C\tau\mE\Bigl[\int_D\int_{t_{n}}^{t_{n+1}}\Bigl|\int_{s}^{t_{n+1}}DF(u(s))G(u(\xi))\, dW(\xi)\Big|^2\, ds\, dx\Bigr]\\\nonumber
		&\qquad+ \frac{1}{8}\mE\bigl[\|e^{n+1} - e^n\|^2_{L^2}\bigr]\\\nonumber
		&=C\tau\mE\Bigl[\int_{t_{n}}^{t_{n+1}}\Bigl\|\int_{s}^{t_{n+1}}DF(u(s))G(u(\xi))\, dW(\xi)\Bigr\|^2_{L^2}\, ds\Bigr] + \frac{1}{8}\mE\bigl[\|e^{n+1} - e^n\|^2_{L^2}\bigr]\\\nonumber
		&=C\tau \mE\Bigl[\int_{t_{n}}^{t_{n+1}} \int_{s}^{t_{n+1}}\|DF(u(s))G(u(\xi))\|^2_{L^2}\, d\xi\, ds\Bigr]+ \frac{1}{8}\mE\bigl[\|e^{n+1} - e^n\|^2_{L^2}\bigr]\\\nonumber
		&\leq C\tau^3 \Bigl(\sup_{s \in [0,T]}\mE\bigl[\|u(s)\|^{4(q-1)}_{L^{4(q-1)}}\bigr]\Bigr)^{1/2}\Bigl(\sup_{\xi\in [0,T]}\mE\bigl[\|u(\xi)\|^4_{H^1}\bigr]\Bigr)^{1/2}\\\nonumber
		&\qquad+ \frac{1}{8}\mE\bigl[\|e^{n+1} - e^n\|^2_{L^2}\bigr]\\\nonumber
		&\leq C\tau^3 + \frac{1}{8}\mE\bigl[\|e^{n+1} - e^n\|^2_{L^2}\bigr].
	\end{align*}
	
	To estimate ${\tt III}$, we use the one-sided Lipschitz condition \eqref{oneside_Lip} as follows.
	\begin{align}
		\mE[{\tt III}]	&\leq C\tau\mE\bigl[\|e^{n+1}\|^2_{L^2}\bigr]\\\nonumber
		&\leq C\tau\mE\bigl[\|e^{n+1} - e^n\|^2_{L^2}\bigr] + C\tau\mE\bigl[\|e^n\|^2_{L^2}\bigr].
	\end{align}
	
	To estimate {\tt IV}, using Lemma \ref{lemma3.1}, the It\^o isometry and the martingale property of It\^o integrals we have
	\begin{align}
		\mE[{\tt IV}] &= \mE\Bigl[\int_{t_{n}}^{t_{n+1}} \bigl(G(u(s)) - \mathcal{G}(s;u^n),e^{n+1} - e^n\bigr)\, dW(s)\Bigr] \\\nonumber
		&\qquad+ \mE\Bigl[\int_{t_{n}}^{t_{n+1}} \bigl(G(u(s)) - \mathcal{G}(s;u^n),e^{n}\bigr)\, dW(s)\Bigr]\\\nonumber
		&= \mE\Bigl[\int_{t_{n}}^{t_{n+1}} \bigl(G(u(s)) - \mathcal{G}(s;u^n),e^{n+1} - e^n\bigr)\, dW(s)\Bigr] + 0 \\\nonumber
		&= \mE\Bigl[\int_{t_{n}}^{t_{n+1}} \bigl(G(u(s)) - \mathcal{G}(s;u(t_n)),e^{n+1} - e^n\bigr)\, dW(s)\Bigr] \\\nonumber
		&\qquad+\mE\Bigl[\int_{t_{n}}^{t_{n+1}} \bigl(\mathcal{G}(s;u(t_n)) - \mathcal{G}(s;u^n),e^{n+1} - e^n\bigr)\, dW(s)\Bigr]\\\nonumber
		&\leq C\mE\Bigl[\Bigl\|\int_{t_n}^{t_{n+1}} \bigl(G(u(s)) - \mathcal{G}(s;u(t_n))\bigr)\, dW(s)\Bigr\|^2_{L^2}\Bigr] \\\nonumber
		&\qquad+C\mE\Bigl[\Bigl\|\int_{t_n}^{t_{n+1}} \bigl(\mathcal{G}(s;u(t_n)) - \mathcal{G}(s;u^n)\bigr)\, dW(s)\Bigr\|^2_{L^2}\Bigr] \\\nonumber
		&\qquad+ \frac18\mE\bigl[\|e^{n+1} - e^n\|^2_{L^2}\bigr]\\\nonumber
		&=C\mE\Bigl[\int_{t_n}^{t_{n+1}}\|G(u(s)) - \mathcal{G}(s;u(t_n))\|^2_{L^2}\, ds\Bigr] \\\nonumber
		&\qquad+ C\mE\Bigl[\int_{t_n}^{t_{n+1}}\|\mathcal{G}(s;u(t_n)) - \mathcal{G}(s;u^n)\|^2_{L^2}\, ds\Bigr]\\\nonumber
		&\qquad+ \frac18\mE\bigl[\|e^{n+1} - e^n\|^2_{L^2}\bigr]\\\nonumber
		&\leq C\tau^{1+ 2(1-\eps)} + C\tau\mE\bigl[\|e^n\|^2_{L^2}\bigr] +  \frac18\mE\bigl[\|e^{n+1} - e^n\|^2_{L^2}\bigr].
	\end{align}
	
	Now, we substitute all the estimates from {\tt I, II, III, IV} into \eqref{eq2.14} and use the left side to absorb the like-terms from the right side of the resulting inequality. In summary, we obtain 
	\begin{align}\label{eq2.25}
		&\frac12\mE\bigl[\|e^{n+1}\|^2_{L^2} - \|e^n\|^2_{L^2}\bigr] + \Bigl(\frac18 - C\tau\Bigr)\mE\bigl[\|e^{n+1} - e^n\|^2_{L^2}\bigr] + \frac{\tau}{2}\mE\bigl[\|\nab e^{n+1}\|^2_{L^2}\bigr]\\\nonumber
		&\leq C\tau^{1+2(1-\eps)} + C\tau\mE\bigl[\|e^{n}\|^2_{L^2}\bigr] + C\tau^3.
	\end{align}
	We choose $\tau \leq \tau_0$ ( for $\tau_0$ small enough) such that $\frac{1}{8} - C\tau \geq 0$, so the middle term on the left side of \eqref{eq2.25} is nonnegtive.
	
	Next, applying the summation $\sum_{n=0}^{m}$ for $0\leq m < N$, we obtain
	\begin{align}
		\mE\bigl[\|e^{m+1}\|^2_{L^2}\bigr] + \tau\sum_{n=0}^m \mE\bigl[\|\nab e^{n+1}\|^2_{L^2}\bigr] \leq C\tau^{2(1-\eps)} + C\tau\sum_{n=0}^{m}\mE\bigl[\|e^n\|^2_{L^2}\bigr].
	\end{align}
	
	By using the discrete Gronwall's inequality and taking supremum over all $0\leq m <M$, we arrive at
	\begin{align}\label{eq2.27}
		\sup_{1\leq n \leq N}\mE\bigl[\|e^{n}\|^2_{L^2}\bigr] + \tau\sum_{n=1}^N \mE\bigl[\|\nab e^{n}\|^2_{L^2}\bigr] \leq Ce^{CT}\tau^{2(1-\eps)}.
	\end{align}

	The proof is complete.
	
\end{proof}

\section{Fully discrete finite element discretization}\label{fully_discrete}
In this section, we consider the $\mathcal{P}_1$-Lagrangian finite element space
\begin{align}\label{eq20180713_1}
	V_h = \bigl\{v_h \in H_0^1(\D): v_h|_{K} \in \mathcal{P}_1(K)\quad\forall K\in\mathcal{T}_h\bigr\},
\end{align}
where $\mathcal{P}_1$ denotes the space of all linear polynomials. Then the finite element approximation of Algorithm 1 is presented in Algorithm 2 as below.

\smallskip
\noindent
\textbf{Algorithm 2} 

We seek an $\mathcal{F}_{t_n}$ adapted $V_h$-valued process $\{u_h^n\}_{n=1}^N$ such that it holds $\mathbb{P}$-almost surely that
\begin{align}\label{dfem}
	&(u^{n+1}_h-u^n_h, v_h) + \tau ( \nabla u^{n+1}_h, \nabla v_h )= \tau (I_hF^{n+1}, v_h)\\
	&\quad+ (G(u^n_h), v_h) \, \Delta W_{n}+ \frac12 DG(u_h^n)\,G(u_h^n)\bigl[(\Delta W_n)^2 - \tau\bigr],v_h\bigr)\qquad \forall \, v_h \in V_h,\notag
\end{align}
where $F^{n+1}:=u^{n+1}_h-(u^{n+1}_h)^q$, ${\Delta} W_{n}=W(t_{n+1})-W(t_n) \sim \mathcal{N} (0,\tau)$, and  $I_h$ is the standard nodal value interpolation operator
$I_h: C(\bar{\D}) \longrightarrow V_h$, i.e.,
\begin{equation} \label{interpolation}
	I_h v := \sum_{i=1}^{N_h} v(a_i)\varphi_i,
\end{equation}
where $N_h$ denotes the number of vertices of the triangulation $\mathcal{T}_h$,
and ${\varphi_i}$ denotes the nodal basis function of $V_h$ corresponding to the vertex $a_i$.
The initial condition is chosen by $u_h^0  = P_h u_0$ where $P_h: L^2(\D) \longrightarrow V_h$ is the $L^2$-projection operator defined by
\begin{align*}
	\bigl(P_h w, v_h\bigr) = (w, v_h) \qquad v_h \in V_h.
\end{align*}
\par
For each $w \in H^s(\D)$  for $s>\frac32$, the following error estimates about the $L^2$-projection can be found in \cite{BS2008,ciarlet2002finite}:
\begin{align}
	\label{Ph1}
	&\|w - P_h w \|_{L^2} + h \| \nabla (w - P_h w) \|_{L^2}
	\leq C h^{\min\{2,s\}} \|w\|_{H^s},\\
	\label{Ph2}
	&\|w - P_h w\|_{L^\infty} \leq C h^{2-\frac{d}{2}} \|w\|_{H^2}.
\end{align}

Finally, given $v_h \in {V}_h$, the discrete Laplace operator $\Delta_h: {V}_h\longrightarrow
{V}_h$ is defined by
\begin{equation} \label{eq:discrete-Laplace}
	(\Delta_h v_h, w_h)=-(\nabla v_h,\nabla w_h) \qquad \forall\, w_h\in V_h.
\end{equation}

\subsection{Stability estimates for the $p$-th moment of the $H^1$-seminorm of $u_h^n$} \label{sec-3.2}
First, we shall prove the second moment discrete $H^1$-seminorm stability result, which is necessary to establish the corresponding higher moment stability result.
\begin{theorem}\label{thm20180711_1}
	Under the mesh constraint \eqref{eq20180907}, we have
	\begin{align}\label{eq20180711_7}
		\sup_{0\leq n \leq N}\E\left[\|\nabla u^{n}_h\|_{L^2}^2\right]&+\frac14\sum_{n=0}^{N-1}\E\left[\|\nabla (u^{n+1}_h-u^{n}_h)\|_{L^2}^2\right] \\
		&\quad +\tau\sum_{n=0}^{N-1}\E\left[\|\Delta_h u_h^{n}\|_{L^2}^2\right]\le C. \notag
	\end{align}
\end{theorem}
\begin{proof}
	Testing \eqref{dfem} with $-\Delta_h u_h^{n+1}$. Then
	\begin{align}\label{eq20180711_1}
		&(u^{n+1}_h-u^n_h, -\Delta_h u_h^{n+1}) + \tau ( \nabla u^{n+1}_h, -\nabla\Delta_h u_h^{n+1} ) \\\nonumber
		&\quad= \tau(I_hF^{n+1}, -\Delta_h u_h^{n+1}) + ( G(u^n_h), -\Delta_h u_h^{n+1}) \, \Delta W_{n+1} \\\nonumber
		&\qquad\qquad+\bigl(\frac12 DG(u_h^n)\,G(u_h^n)((\Delta W_n)^2 - \tau),-\Delta_h u_h^{n+1}\bigr).\notag
	\end{align}
	
	Using the definition of the discrete Laplace operator and the simple identity $2a(a-b) = a^2 - b^2 + (a-b)^2$, we get
	\begin{align}\label{eq20180711_2}
		(u^{n+1}_h-u^n_h, -\Delta_h u_h^{n+1})&=\frac12\|\nabla u^{n+1}_h\|_{L^2}^2-\frac12\|\nabla u^{n}_h\|_{L^2}^2\\
		&\qquad+\frac12\|\nabla (u^{n+1}_h-u^{n}_h)\|_{L^2}^2,\notag\\
		\tau ( \nabla u^{n+1}_h, -\nabla\Delta_h u_h^{n+1} )&=\tau\|\Delta_h u_h^{n+1}\|_{L^2}^2.\label{eq20181009_2}
	\end{align}
	
	The expectation of the second term on the right-hand side of \eqref{eq20180711_1} can be bounded by
	\begin{align}
		\E[( G(u^n_h), -\Delta_h u_h^{n+1}) \, \Delta W_{n}]&=\E [(\nabla(P_hG(u^n_h)), \nabla (u_h^{n+1}-u^n_h)) \, {\Delta} W_{n}]\label{eq20181009_3}\\
		&\le C\tau\E[\|\nabla u^n_h\|_{L^2}^2]+\frac14\E[\|\nabla (u_h^{n+1}-u^n_h)\|_{L^2}^2].\notag
	\end{align}
	
	The expectation of the third term on the right-hand side of \eqref{eq20180711_1} can be bounded by
	\begin{align}\label{eq20221026_1}
		&\frac12\E[( DG(u_h^n)\,G(u_h^n)((\Delta W_n)^2 - \tau),-\Delta_h u_h^{n+1})]\\
		=&\frac12\E [(\nabla(P_h(DG(u_h^n)\,G(u_h^n))), \nabla (u_h^{n+1}-u^n_h)) \, ((\Delta W_n)^2 - \tau)]\notag\\
		\le&C\tau^2\E[\|\nabla u^n_h\|_{L^2}^2]+\frac14\E[\|\nabla (u_h^{n+1}-u^n_h)\|_{L^2}^2],\notag
	\end{align}
	where the last inequality is obtained by using the assumption {\bf(A2)}. Notice that the stability in the $H^1$-seminorm of the $L^2$-projection (see \cite{bank2014h}) is used in the inequalities of \eqref{eq20181009_3} and \eqref{eq20221026_1}.
	
	For the first term on the right-hand side of \eqref{eq20180711_1} since it cannot be treated as a bad term, which aligns with the continuous case. Denote $u_i=u_h^{n+1}(a_i)$, and then
	\begin{align}\label{eq20180711_3}
		\tau (I_hF^{n+1}, -\Delta_h u_h^{n+1})&=\tau\|\nabla u^{n+1}_h\|_{L^2}^2-\tau(\nabla\sum_{i=1}^{N_h} u_i^q \varphi_i,\nabla \sum_{j=1}^{N_h} u_j\varphi_j)\\
		&=\tau\|\nabla u^{n+1}_h\|_{L^2}^2-\tau \sum_{i,j=1}^{N_h} ( u_i^q \nabla\varphi_i,  u_j\nabla\varphi_j)\notag\\
		&=\tau\|\nabla u^{n+1}_h\|_{L^2}^2-\tau \sum_{i,j=1}^{N_h} b_{ij}(\nabla\varphi_i,  \nabla\varphi_j)\notag,
	\end{align}
	where $b_{ij}=u_i^q u_j$.
	
	Using Young's inequality when $i\neq j$, we have
	\begin{align}\label{eq20180711_4}
		|b_{ij}|\le \frac{q}{q+1}u_i^{q+1}+\frac{1}{q+1}u_j^{q+1}.
	\end{align}
	
	Besides, since the stiffness matrix is diagonally dominant, we have
	\begin{align}
		-\tau \sum_{i,j=1}^{N_h} b_{ij}(\nabla\varphi_i,  \nabla\varphi_j)&\le-\tau \sum_{k=1}^{N_h} b_{kk}[(\nabla\varphi_k,  \nabla\varphi_k)-\frac{q}{q+1}\sum_{i=1,\atop i\neq k}^{N_h} |(\nabla\varphi_i,  \nabla\varphi_k)|\\
		&\quad-\frac{1}{q+1}\sum_{j=1,\atop j\neq k}^{N_h} |(\nabla\varphi_k,  \nabla\varphi_j)|]\notag\\
		&\le-\tau \sum_{k=1}^{N_h} b_{kk}[(\nabla\varphi_k,  \nabla\varphi_k)-\sum_{i=1,\atop i\neq k}^{N_h} (\nabla\varphi_i,  \nabla\varphi_k)]\notag\\
		&\le0\notag.
	\end{align}

	Then we have
	\begin{align}\label{eq20180711_6}
		\tau (I_hF^{n+1}, -\Delta_h u_h^{n+1})\le\tau\|\nabla u^{n+1}_h\|_{L^2}^2.
	\end{align}
	
	Combining \eqref{eq20180711_1}--\eqref{eq20181009_3} and \eqref{eq20180711_6}, and taking the summation, we have
	\begin{align}\label{eq20180712_1}
		&\frac12\E\left[\|\nabla u^{\ell}_h\|_{L^2}^2\right]+\frac14\sum_{n=0}^{\ell-1}\E\left[\|\nabla (u^{n+1}_h-u^{n}_h)\|_{L^2}^2\right]+\tau\sum_{n=0}^{\ell-1}\E\left[\|\Delta_h u_h^{n+1}\|_{L^2}^2\right]\\
		&\quad\le C\tau\sum_{n=0}^{\ell-1}\E[\|\nabla u^n_h\|_{L^2}^2].\notag
	\end{align}
	
	Using Gronwall's inequality, we obtain \eqref{eq20180711_7}.
\end{proof}

Before we establish the error estimates, we need to prove the stability of
the higher moments for the $H^1$-seminorm of the numerical solution.

\begin{theorem}\label{thm20180802_1}
	Suppose the mesh assumption \eqref{eq20180907} holds. Then for any $p\ge2$,
	\begin{align*}
		\sup_{0 \leq n \leq M} \E\left[\|\nabla u^{n}_h\|_{L^2}^p\right]\le C.
	\end{align*}
\end{theorem}
\begin{proof}
	The proof is divided into three steps. In Step 1, we establish the bound for $\E\|\nabla u^{\ell}_h\|_{L^2}^{4}$.
	In Step 2, we give the bound for $\E\|\nabla u^{\ell}_h\|_{L^2}^p$, where $p=2^r$ and $r$ is an arbitrary
	positive integer. In Step 3, we obtain the bound for $\E\|\nabla u^{\ell}_h\|_{L^2}^p$, where $p$ is an
	arbitrary real number and $p\ge2$.
	
	\smallskip
	{\bf Step 1.} Based on \eqref{eq20180711_1}--\eqref{eq20180711_6}, we have
	\begin{align}\label{eq20180802_1}
		&\frac12 \|\nabla u^{n+1}_h\|_{L^2}^2 -\frac12 \|\nabla u^{n}_h\|_{L^2}^2 +\frac12 \|\nabla (u^{n+1}_h-u^{n}_h)\|_{L^2}^2 +\tau \|\Delta_h u_h^{n+1}\|_{L^2}^2 \\[2mm]
		& -( G(u^n_h), -\Delta_h u_h^{n+1}) \, \Delta W_{n} -\frac12\bigl(DG(u_h^n)\,G(u_h^n)((\Delta W_n)^2 - \tau),-\Delta_h u_h^{n+1}\bigr)\notag\\
		&\le\tau \|\nabla u^{n+1}_h\|_{L^2}^2.\notag
	\end{align}
	
	Note the following identity
	\begin{align}\label{eq20180802_2}
		\|\nabla u^{n+1}_h\|_{L^2}^2+\frac12\|\nabla u^{n}_h\|_{L^2}^2&=\frac34(\|\nabla u^{n+1}_h\|_{L^2}^2+\|\nabla u^n_h\|_{L^2}^2)\\
		&\quad+\frac14(\|\nabla u^{n+1}_h\|_{L^2}^2-\|\nabla u^n_h\|_{L^2}^2).\notag
	\end{align}
	Multiplying \eqref{eq20180802_1} by $\|\nabla u^{n+1}_h\|_{L^2}^2+\frac12\|\nabla u^{n}_h\|_{L^2}^2$, we obtain
	\begin{align}\label{eq20180802_3}
		&\frac38(\|\nabla u^{n+1}_h\|_{L^2}^4-\|\nabla u^n_h\|_{L^2}^4)+\frac18(\|\nabla u^{n+1}_h\|_{L^2}^2-\|\nabla u^n_h\|_{L^2}^2)^2\\
		&\quad+(\frac12\|\nabla (u^{n+1}_h-u^{n}_h)\|_{L^2}^2+\tau\|\Delta_h u_h^{n+1}\|_{L^2}^2)(\|\nabla u^{n+1}_h\|_{L^2}^2+\frac12\|\nabla u^{n}_h\|_{L^2}^2)\notag\\
		&\le\tau\|\nabla u^{n+1}_h\|_{L^2}^2(\|\nabla u^{n+1}_h\|_{L^2}^2+\frac12\|\nabla u^{n}_h\|_{L^2}^2)\notag\\
		&\quad+( G(u^n_h), -\Delta_h u_h^{n+1}) \, {\Delta} W_{n}(\|\nabla u^{n+1}_h\|_{L^2}^2+\frac12\|\nabla u^{n}_h\|_{L^2}^2)\notag\\
		&\quad+\frac12\bigl(DG(u_h^n)\,G(u_h^n)((\Delta W_n)^2 - \tau),-\Delta_h u_h^{n+1}\bigr)(\|\nabla u^{n+1}_h\|_{L^2}^2+\frac12\|\nabla u^{n}_h\|_{L^2}^2).\notag
	\end{align}
	
	The first term on the right-hand side of \eqref{eq20180802_3} can be written as
	\begin{align}\label{eq20180802_4}
		&\tau\|\nabla u^{n+1}_h\|_{L^2}^2(\|\nabla u^{n+1}_h\|_{L^2}^2+\frac12\|\nabla u^{n}_h\|_{L^2}^2)\\
		&\quad=\tau\|\nabla u^{n+1}_h\|_{L^2}^2(\frac32\|\nabla u^{n+1}_h\|_{L^2}^2-\frac12(\|\nabla u^{n+1}_h\|_{L^2}^2-\|\nabla u^n_h\|_{L^2}^2))\notag\\
		&\quad\le C\tau\|\nabla u^{n+1}_h\|_{L^2}^4+\theta_1(\|\nabla u^{n+1}_h\|_{L^2}^2-\|\nabla u^n_h\|_{L^2}^2)^2\notag,
	\end{align}
	where $\theta_1>0$ will be determined later.
	
	The second term on the right-hand side of \eqref{eq20180802_3} can be written as
	\begin{align}\label{eq20180802_5}
		&(G(u^n_h), -\Delta_h u_h^{n+1}) \, {\Delta} W_{n}(\|\nabla u^{n+1}_h\|_{L^2}^2+\frac12\|\nabla u^{n}_h\|_{L^2}^2)\\
		&\quad=(\nabla P_hG(u^n_h), \nabla u_h^{n+1}) \, {\Delta} W_{n}(\|\nabla u^{n+1}_h\|_{L^2}^2+\frac12\|\nabla u^{n}_h\|_{L^2}^2)\notag\\
		&\quad= ((\nabla P_hG(u^n_h), \nabla u_h^{n+1}-\nabla u_h^n){\Delta} W_{n}\notag\\
		&\qquad+(\nabla P_hG(u^n_h),\nabla u_h^n){\Delta} W_{n})(\|\nabla u^{n+1}_h\|_{L^2}^2+\frac12\|\nabla u^{n}_h\|_{L^2}^2)\notag\\
		&\quad\le(\frac14\|\nabla u_h^{n+1}-\nabla u_h^n\|_{L^2}^2+C\|\nabla u_h^n\|_{L^2}^2({\Delta} W_{n})^2\notag\\
		&\qquad+(\nabla P_hG(u^n_h),\nabla u_h^n){\Delta} W_{n})(\|\nabla u^{n+1}_h\|_{L^2}^2+\frac12\|\nabla u^{n}_h\|_{L^2}^2)\notag.
	\end{align}
	
	For the right-hand side of \eqref{eq20180802_5}, using the Cauchy-Schwarz inequality, we get
	\begin{align}\label{eq20180802_6}
		&C\|\nabla u_h^n\|_{L^2}^2({\Delta} W_{n})^2(\|\nabla u^{n+1}_h\|_{L^2}^2+\frac12\|\nabla u^{n}_h\|_{L^2}^2)\\
		&\quad=C\|\nabla u_h^n\|_{L^2}^2({\Delta} W_{n})^2(\|\nabla u^{n+1}_h\|_{L^2}^2-\|\nabla u^{n}_h\|_{L^2}^2+\frac32\|\nabla u^{n}_h\|_{L^2}^2)\notag\\
		&\quad\le\theta_2(\|\nabla u^{n+1}_h\|_{L^2}^2-\|\nabla u^n_h\|_{L^2}^2)^2+C\|\nabla u_h^n\|_{L^2}^4({\Delta} W_{n})^4\notag\\
		&\qquad+C\|\nabla u_h^n\|_{L^2}^4({\Delta} W_{n})^2,\notag
	\end{align}
	where $\theta_2>0$ will be determined later.
	
	Similarly, using the Cauchy-Schwarz inequality, we have
	\begin{align}\label{eq20180802_7}
		&(\nabla P_hG(u^n_h),\nabla u_h^n){\Delta} W_{n}(\|\nabla u^{n+1}_h\|_{L^2}^2+\frac12\|\nabla u^{n}_h\|_{L^2}^2)\\
		&\quad=(\nabla P_hG(u^n_h),\nabla u_h^n){\Delta} W_{n}(\|\nabla u^{n+1}_h\|_{L^2}^2-\|\nabla u^{n}_h\|_{L^2}^2+\frac32\|\nabla u^{n}_h\|_{L^2}^2)\notag\\
		&\quad\le\theta_3(\|\nabla u^{n+1}_h\|_{L^2}^2-\|\nabla u^n_h\|_{L^2}^2)^2+C\|\nabla u_h^n\|_{L^2}^4({\Delta} W_{n})^2\notag\\
		&\qquad+\frac32(\nabla P_hG(u^n_h),\nabla u_h^n){\Delta} W_{n}\|\nabla u^{n}_h\|_{L^2}^2\notag,
	\end{align}
	where $\theta_3>0$ will be determined later.
	
	The third term on the right-hand side of \eqref{eq20180802_3} can be written as
	\begin{align}\label{eq20221226_2}
		&\quad\frac12\bigl(DG(u_h^n)\,G(u_h^n)((\Delta W_n)^2 - \tau),-\Delta_h u_h^{n+1}\bigr)(\|\nabla u^{n+1}_h\|_{L^2}^2+\frac12\|\nabla u^{n}_h\|_{L^2}^2)\\
		&=\frac12\bigl(\nabla P_h(DG(u_h^n)\,G(u_h^n))((\Delta W_n)^2 - \tau),\nabla(u_h^{n+1}-u_h^n)\bigr)\notag\\
		&\quad(\|\nabla u^{n+1}_h\|_{L^2}^2+\frac12\|\nabla u^{n}_h\|_{L^2}^2) +\frac12\bigl(\nabla P_h(DG(u_h^n)\,G(u_h^n))((\Delta W_n)^2 - \tau),\nabla u_h^n\bigr)\notag\\
		&\quad(\|\nabla u^{n+1}_h\|_{L^2}^2+\frac12\|\nabla u^{n}_h\|_{L^2}^2)\notag\\
		&\le(\frac14\|\nabla u_h^{n+1}-\nabla u_h^n\|_{L^2}^2+C\|\nabla u_h^n\|_{L^2}^2((\Delta W_n)^2-\tau)^2\notag\\
		&\quad +\frac12\bigl(\nabla P_h(DG(u_h^n)\,G(u_h^n))((\Delta W_n)^2 - \tau),\nabla u_h^n\bigr)(\|\nabla u^{n+1}_h\|_{L^2}^2+\frac12\|\nabla u^{n}_h\|_{L^2}^2)\notag.
	\end{align}
	
	For the right-hand side of \eqref{eq20221226_2}, using the Cauchy-Schwarz inequality, we get
	\begin{align}
		&C\|\nabla u_h^n\|_{L^2}^2((\Delta W_n)^2-\tau)^2(\|\nabla u^{n+1}_h\|_{L^2}^2+\frac12\|\nabla u^{n}_h\|_{L^2}^2)\\
		&\quad=C\|\nabla u_h^n\|_{L^2}^2((\Delta W_n)^2-\tau)^2(\|\nabla u^{n+1}_h\|_{L^2}^2-\|\nabla u^{n}_h\|_{L^2}^2+\frac32\|\nabla u^{n}_h\|_{L^2}^2)\notag\\
		&\quad\le\theta_4(\|\nabla u^{n+1}_h\|_{L^2}^2-\|\nabla u^n_h\|_{L^2}^2)^2+C\|\nabla u_h^n\|_{L^2}^4((\Delta W_n)^2-\tau)^4 \notag\\
		&\qquad +C\|\nabla u_h^n\|_{L^2}^4((\Delta W_n)^2-\tau)^2,\notag
	\end{align}
	where $\theta_4>0$ will be determined later.
	Similarly, using the Cauchy-Schwarz inequality, we have
	\begin{align}
		&(\nabla P_h(DG(u_h^n)\,G(u_h^n)),\nabla u_h^n)((\Delta W_n)^2 - \tau)(\|\nabla u^{n+1}_h\|_{L^2}^2+\frac12\|\nabla u^{n}_h\|_{L^2}^2)\\
		&\quad=(\nabla P_h(DG(u_h^n)\,G(u_h^n)),\nabla u_h^n)((\Delta W_n)^2 - \tau)(\|\nabla u^{n+1}_h\|_{L^2}^2\notag\\
		&\qquad-\|\nabla u^{n}_h\|_{L^2}^2+\frac32\|\nabla u^{n}_h\|_{L^2}^2)\notag\\
		&\quad\le\theta_5(\|\nabla u^{n+1}_h\|_{L^2}^2-\|\nabla u^n_h\|_{L^2}^2)^2+C\|\nabla u_h^n\|_{L^2}^4((\Delta W_n)^2 - \tau)^2\notag\\
		&\qquad+\frac32(\nabla P_h(DG(u^n)\,G(u^n)),\nabla u_h^n)((\Delta W_n)^2 - \tau)\|\nabla u^{n}_h\|_{L^2}^2\notag,
	\end{align}
	where $\theta_5>0$ will be determined later.
	
	Choosing $\theta_1\sim \theta_5$ such that $\theta_1+\cdots+\theta_3\le\frac{1}{16}$, and then taking the summation over $n$ from $0$ to $\ell-1$ and taking the expectation on both sides of \eqref{eq20180802_3}, we obtain
	\begin{align}\label{eq20180802_8}
		&\frac38\E\left[\|\nabla u^{\ell}_h\|_{L^2}^4\right]+\frac{1}{16}\sum_{n=0}^{\ell-1}\E\left[(\|\nabla u^{n+1}_h\|_{L^2}^2-\|\nabla u^n_h\|_{L^2}^2)^2\right]\\
		&+\sum_{n=0}^{\ell-1}\E\left[(\frac14\|\nabla (u^{n+1}_h-u^{n}_h)\|_{L^2}^2+\tau\|\Delta_h u_h^{n+1}\|_{L^2}^2)(\|\nabla u^{n+1}_h\|_{L^2}^2+\frac12\|\nabla u^{n}_h\|_{L^2}^2)\right]\notag\\
		&\le C\tau\sum_{n=0}^{\ell-1}\E\left[\|\nabla u^{n+1}_h\|_{L^2}^4\right]+\frac38\E\left[\|\nabla u^0_h\|_{L^2}^4\right]+C\tau^2\sum_{n=0}^{\ell-1}\E\left[\|\nabla u_h^n\|_{L^2}^4\right]\notag\\
		&+C\tau\sum_{n=0}^{\ell-1}\E\left[\|\nabla u_h^n\|_{L^2}^4\right].\notag
	\end{align}
	
	When restricting $\tau\le C$, we have
	\begin{align}\label{eq20180802_9}
		&\frac14\E\left[\|\nabla u^{\ell}_h\|_{L^2}^4\right]+\frac{1}{16}\sum_{n=0}^{\ell-1}\E\left[(\|\nabla u^{n+1}_h\|_{L^2}^2-\|\nabla u^n_h\|_{L^2}^2)^2\right]\\
		&+\sum_{n=0}^{\ell-1}\E\left[(\frac14\|\nabla (u^{n+1}_h-u^{n}_h)\|_{L^2}^2+\tau\|\Delta_h u_h^{n+1}\|_{L^2}^2)(\|\nabla u^{n+1}_h\|_{L^2}^2+\frac12\|\nabla u^{n}_h\|_{L^2}^2)\right]\notag\\
		&\le C\tau\sum_{n=0}^{\ell-1}\E\left[\|\nabla u^{n}_h\|_{L^2}^4\right]+\frac38\E\left[\|\nabla u^0_h\|_{L^2}^4\right].\notag
	\end{align}
	
	Using Gronwall's inequality, we obtain
	\begin{align}\label{eq20180802_10}
		&\frac14\E\left[\|\nabla u^{\ell}_h\|_{L^2}^4\right]+\frac{1}{16}\sum_{n=0}^{\ell-1}\E\left[(\|\nabla u^{n+1}_h\|_{L^2}^2-\|\nabla u^n_h\|_{L^2}^2)^2\right]\\
		&\qquad+\sum_{n=0}^{\ell-1}\E\bigl[(\frac14\|\nabla (u^{n+1}_h-u^{n}_h)\|_{L^2}^2+\tau\|\Delta_h u_h^{n+1}\|_{L^2}^2)(\|\nabla u^{n+1}_h\|_{L^2}^2\notag\\
		&\qquad+\frac12\|\nabla u^{n}_h\|_{L^2}^2)\bigr]\le C.\notag
	\end{align}
	
	\smallskip
	{\bf Step 2.} Similar to Step 1, using \eqref{eq20180802_3}--\eqref{eq20180802_7}, we have
	\begin{align}\label{eq20180802_11}
		&\frac38(\|\nabla u^{n+1}_h\|_{L^2}^4-\|\nabla u^n_h\|_{L^2}^4)+\frac{1}{16}(\|\nabla u^{n+1}_h\|_{L^2}^2-\|\nabla u^n_h\|_{L^2}^2)^2\\
		&\quad+(\frac14\|\nabla (u^{n+1}_h-u^{n}_h)\|_{L^2}^2+\tau\|\Delta_h u_h^{n+1}\|_{L^2}^2)(\|\nabla u^{n+1}_h\|_{L^2}^2+\frac12\|\nabla u^{n}_h\|_{L^2}^2)\notag\\
		&\le C\tau\|\nabla u^{n+1}_h\|_{L^2}^4+C\|\nabla u_h^n\|_{L^2}^4({\Delta} W_{n})^4+C\|\nabla u_h^n\|_{L^2}^4({\Delta} W_{n})^2\notag\\
		&\quad+C\|\nabla u_h^n\|_{L^2}^4{\Delta} W_{n}.\notag
	\end{align}
	
	Proceed similarly as in Step 1. Multiplying \eqref{eq20180802_11} with $\|\nabla u^{n+1}_h\|_{L^2}^4+\frac12\|\nabla u^{n}_h\|_{L^2}^4$, we can obtain the 8-th moment of the $H^1$-seminorm stability result of the numerical solution.
	Then repeat this process, the $2^r$-th moment of the $H^1$-seminorm stability result of the numerical solution
	can be obtained.
	
	\smallskip
	{\bf Step 3.} Suppose $2^{r-1}\le p\le 2^r$. By Young's inequality, we have
	\begin{align}\label{eq20180802_12}
		\E\left[\|\nabla u^{\ell}_h\|_{L^2}^p\right]&\le \E\left[\|\nabla u^{\ell}_h\|_{L^2}^{2^r}\right]+C< \infty,
	\end{align}
	where the second inequality follows from the results of Step 2. The proof is completed.
\end{proof}

\subsection{Stability estimates for the $p$-th moment of the $L^2$-norm of $u_h^n$}
Since the mass matrix may not be the diagonally dominated matrix, we cannot use the above idea to prove the $L^2$ stability. Instead, we prove the stability results by utilizing the above established results. The following results hold when $q\ge3$ is the odd integer in 2D case, and when $q=3$ or $q=5$ in 3D case.
\begin{theorem}\label{thm20180911}
	Under the mesh assumption \eqref{eq20180907},   there holds
	\begin{align*}
		\sup_{0\leq n \leq N}\E\left[\|u^{n}_h\|_{L^2}^2\right]+\sum_{n=0}^{N-1}\E\left[\|(u^{n+1}_h-u^{n}_h)\|_{L^2}^2\right] &+\tau\sum_{n=0}^{N-1}\E\left[\|\nabla u_h^{n+1}\|_{L^2}^2\right]\\
		&+\frac{\tau}{2}\sum_{n=0}^{N-1}\E\left[\|u_h^{n+1}\|_{L^{q+1}}^{q+1}\right]\le C.
	\end{align*}
\end{theorem}
\begin{proof}
	Testing \eqref{dfem} with $u_h^{n+1}$ yields
	\begin{align}\label{eq20221227_1}
		&(u^{n+1}_h-u^n_h, u_h^{n+1}) + \tau ( \nabla u^{n+1}_h, \nabla u_h^{n+1} )= \tau (I_hF^{n+1}, u_h^{n+1})\\
		&\quad+ (G(u^n_h), u_h^{n+1}) \, \Delta W_{n} + \frac12 DG(u_h^n)\,G(u_h^n)\bigl[(\Delta W_n)^2 - \tau\bigr],u_h^{n+1}\bigr).\notag
	\end{align}
	
	We can easily prove the following inequalities:
	\begin{align*}
		(u^{n+1}_h-u^n_h, u_h^{n+1})&=\frac12\| u^{n+1}_h\|_{L^2}^2-\frac12\| u^{n}_h\|_{L^2}^2+\frac12\|u^{n+1}_h-u^{n}_h\|_{L^2}^2,\\
		\E[(G(u^n_h),  u_h^{n+1}) \, \Delta W_{n}]&=\E [(G(u^n_h),  (u_h^{n+1}-u^n_h)) \, \Delta W_{n}]\\
		&\le C\tau+C\tau\E[\| u^n_h\|_{L^2}^2]+\frac14\E[\|u_h^{n+1}-u^n_h\|_{L^2}^2],\\
		\E[DG(u_h^n)\,G(u_h^n)((\Delta W_n)^2 - \tau),u_h^{n+1}\bigr)]&=\E[DG(u_h^n)\,G(u_h^n)((\Delta W_n)^2 - \tau),u_h^{n+1}-u_h^n\bigr)]\\
		&\le C\tau^2+C\tau^2\E[\| u^n_h\|_{L^2}^2]+\frac14\E[\|u_h^{n+1}-u^n_h\|_{L^2}^2],
	\end{align*}
	where {\bf (A2)} is used in the inequality above.
	
	We have the following standard interpolation result and the inverse inequality (see \cite{ciarlet2002finite}):
	\begin{align}\label{eq20180807_7}
		\|v-I_hv\|_{L^{\frac{q+1}{q}}(K)} &\le Ch_K\|\nabla v\|_{L^{\frac{q+1}{q}}(K)},\\
		\|v\|_{L^{q+1}(K)}^{q+1} &\le\frac{C}{h_K^{d\cdot\frac{q-1}{2}}}\|v\|_{L^2(K)}^{q+1}.\label{eq20180807_8}
	\end{align}
	
	Using \eqref{eq20180807_7}--\eqref{eq20180807_8}, and Young's inequality, we have
	\begin{align}\label{eq20180808_2}
		&\tau(I_hF^{n+1}, u_h^{n+1})=\tau(F^{n+1}, u_h^{n+1})-\tau(F^{n+1}-I_hF^{n+1}, u_h^{n+1})\\
		&\quad\le\tau\|u_h^{n+1}\|_{L^2}^2-\tau\|u_h^{n+1}\|_{L^{q+1}}^{q+1}\notag\\
		&\qquad+C\tau\|F^{n+1}-I_hF^{n+1}\|_{L^{\frac{q+1}{q}}}^{\frac{q+1}{q}}+\frac{\tau}{4}\|u_h^{n+1}\|_{L^{q+1}}^{q+1}\notag\\
		&\quad\le\tau\|u_h^{n+1}\|_{L^2}^2-\tau\|u_h^{n+1}\|_{L^{q+1}}^{q+1}\notag\\
		&\qquad+C\tau \sum_{K\in\mathcal{T}_h}h_K^{\frac{q+1}{q}}\bigl((u_h^{n+1})^{\frac{q^2-1}{q}},(\nabla u_h^{n+1})^{\frac{q+1}{q}}\bigr)_K+\frac{\tau}{4}\|u_h^{n+1}\|_{L^{q+1}}^{q+1}\notag\\
		&\quad\le\tau\|u_h^{n+1}\|_{L^2}^2-\frac{\tau}{2}\|u_h^{n+1}\|_{L^{q+1}}^{q+1}+C\tau \sum_{K\in\mathcal{T}_h}h_K^{q+1}\|\nabla u_h^{n+1}\|_{L^{q+1}(K)}^{q+1}\notag\\
		&\quad\le\tau\|u_h^{n+1}\|_{L^2}^2-\frac{\tau}{2}\|u_h^{n+1}\|_{L^{q+1}}^{q+1}+C\tau\sum_{K\in\mathcal{T}_h}h_K^{q+1-d\frac{q-1}{2}}\|\nabla u_h^{n+1}\|_{L^2(K)}^{q+1}.\notag
	\end{align}
	
	Note when $d=2$, $q+1-d\frac{q-1}{2}\ge0$ if $q\ge0$, and when $d=3$, $q+1-d\frac{q-1}{2}\ge0$ if $q\le5$. Using the above inequalities, Theorem \ref{thm20180802_1}, taking summation over $n$ from $0$ to $\ell-1$, and taking expectation on both sides of \eqref{eq20221227_1}, we obtain
	\begin{align}\label{eq20180807_1}
		&\frac14\E\left[\|u^{\ell}_h\|_{L^2}^2\right]+\frac14\sum_{n=0}^{\ell-1}\E\left[\|(u^{n+1}_h-u^{n}_h)\|_{L^2}^2\right]+\tau\sum_{n=0}^{\ell-1}\E\left[\|\nabla u_h^{n+1}\|_{L^2}^2\right]\\
		&\qquad+\frac{\tau}{2}\sum_{n=0}^{\ell-1}\E\left[\|u_h^{n+1}\|_{L^{q+1}}^{q+1}\right]\notag\\
		&\le\tau\sum_{n=0}^{\ell-1}\E\left[\|u^n_h\|_{L^2}^2\right]+C\tau\sum_{n=0}^{\ell-1}\E\left[\|\nabla u_h^{n+1}\|_{L^2}^{q+1}\right]+C\notag\\
		&\le\tau\sum_{n=0}^{\ell-1}\E\left[\|u^n_h\|_{L^2}^2\right]+C,\notag
	\end{align}
	where Theorem \ref{thm20180802_1} is used in the last inequality.
	
	The conclusion is a direct result by using Gronwall's inequality.
\end{proof}

%
%

To obtain the error estimates results, we need to establish a higher moment discrete $L^2$ stability result
for the numerical solution $u_h$.

\begin{theorem}\label{thm20180808_1}
	Suppose the mesh assumption \eqref{eq20180907} holds. Then there holds for any $p\ge2$,
	\begin{align*}
		\sup_{0\leq \ell \leq N}\E\left[\|u^{\ell}_h\|_{L^2}^p\right]\le C\notag.
	\end{align*}
\end{theorem}

\begin{proof}
	The proof is divided into three steps. In Step 1, we give the bound for $\E\| u^{\ell}_h\|_{L^2}^{4}$.
	In Step 2, we give the bound for $\E\| u^{\ell}_h\|_{L^2}^p$, where $p=2^r$ and $r$ is an arbitrary positive integer.
	In Step 3, we give the bound for $\E\| u^{\ell}_h\|_{L^2}^p$, where $p$ is an arbitrary real number and $p\ge2$.
	
	\smallskip
	{\rm Step 1.} Based on \eqref{eq20221227_1}--\eqref{eq20180808_2}, we have
	\begin{align}\label{eq20180808_3}
		\frac12\|u^{n+1}_h\|_{L^2}^2&-\frac12\|u^{n}_h\|_{L^2}^2 +\frac12\|u^{n+1}_h-u^{n}_h\|_{L^2}^2 +\tau\|\nabla u_h^{n+1}\|_{L^2}^2+\frac{\tau}{2}\|u_h^{n+1}\|_{L^{q+1}}^{q+1}\\[2mm]
		& \le\tau\|u_h^{n+1}\|_{L^2}^2+C\tau\|\nabla u_h^{n+1}\|_{L^2}^{q+1}+(G(u^n_h), u_h^{n+1}) \, \Delta W_{n}\notag\\
		&\qquad + \frac12 \bigl(DG(u_h^n)\,G(u_h^n)((\Delta W_n)^2 - \tau),u_h^{n+1}\bigr).\notag
	\end{align}
	
	Note the following identity
	\begin{align}\label{eq20180808_4}
		\|u^{n+1}_h\|_{L^2}^2+\frac12\|u^{n}_h\|_{L^2}^2=&\frac34(\|u^{n+1}_h\|_{L^2}^2+\|u^n_h\|_{L^2}^2)+\frac14(\|u^{n+1}_h\|_{L^2}^2-\|u^n_h\|_{L^2}^2).
	\end{align}
	Multiplying \eqref{eq20180808_3} by $\|u^{n+1}_h\|_{L^2}^2+\frac12\|u^{n}_h\|_{L^2}^2$, we obtain
	\begin{align}\label{eq20180808_5}
		&\frac38(\|u^{n+1}_h\|_{L^2}^4-\|u^n_h\|_{L^2}^4)+\frac18(\|u^{n+1}_h\|_{L^2}^2-\|u^n_h\|_{L^2}^2)^2+(\frac12\|(u^{n+1}_h-u^{n}_h)\|_{L^2}^2\\
		&\quad+\tau\|\nabla u_h^{n+1}\|_{L^2}^2+\frac{\tau}{2}\|u_h^{n+1}\|_{L^{q+1}}^{q+1})(\|u^{n+1}_h\|_{L^2}^2+\frac12\|u^{n}_h\|_{L^2}^2)\notag\\
		&\le(\tau\|u_h^{n+1}\|_{L^2}^2+C\tau\|\nabla u_h^{n+1}\|_{L^2}^{q+1})(\|u^{n+1}_h\|_{L^2}^2+\frac12\|u^{n}_h\|_{L^2}^2)\notag\\
		&\quad+(G(u^n_h), u_h^{n+1})\, \Delta W_{n}(\|u^{n+1}_h\|_{L^2}^2+\frac12\|u^{n}_h\|_{L^2}^2)\notag\\
		&\quad+\frac12 \bigl(DG(u_h^n)\,G(u_h^n)((\Delta W_n)^2 - k),u_h^{n+1}\bigr)(\|u^{n+1}_h\|_{L^2}^2+\frac12\|u^{n}_h\|_{L^2}^2).\notag
	\end{align}
	
	The first term on the right-hand side of \eqref{eq20180808_5} can be written as
	\begin{align}\label{eq20180808_6}
		&(\tau\|u_h^{n+1}\|_{L^2}^2+C\tau\|\nabla u_h^{n+1}\|_{L^2}^{q+1})(\|u^{n+1}_h\|_{L^2}^2+\frac12\|u^{n}_h\|_{L^2}^2)\\
		&\le\tau\|u^{n+1}_h\|_{L^2}^2(\frac32\|u^{n+1}_h\|_{L^2}^2-\frac12(\|u^{n+1}_h\|_{L^2}^2-\|u^n_h\|_{L^2}^2))\notag\\
		&\quad+C\tau\|\nabla u_h^{n+1}\|_{L^2}^{2(q+1)}+\tau\|u^{n+1}_h\|_{L^2}^4+\tau(\|u^{n+1}_h\|_{L^2}^2-\|u^n_h\|_{L^2}^2)^2\notag\\
		&\le C\tau\|u^{n+1}_h\|_{L^2}^4+C\tau\|\nabla u_h^{n+1}\|_{L^2}^{2(q+1)}+\theta_1(\|u^{n+1}_h\|_{L^2}^2-\|u^n_h\|_{L^2}^2)^2\notag,
	\end{align}
	where $\theta_1>0$ will be determined later.
	
	The second term on the right-hand side of \eqref{eq20180808_5} can be written as
	\begin{align}\label{eq20180808_7}
		&(G(u^n_h), u_h^{n+1})\, {\Delta} W_{n}(\|u^{n+1}_h\|_{L^2}^2+\frac12\|u^{n}_h\|_{L^2}^2)\\
		&=(G(u^n_h), u_h^{n+1}-u_h^n+u_h^n) \, {\Delta} W_{n}(\|u^{n+1}_h\|_{L^2}^2+\frac12\|u^{n}_h\|_{L^2}^2)\notag\\
		&\le(\frac14\|u_h^{n+1}-u_h^n\|_{L^2}^2+C(1+\|u_h^n\|_{L^2}^2)({\Delta} W_{n})^2\notag\\
		&\quad+(G(u_h^n),u_h^n) {\Delta} W_{n})(\|u^{n+1}_h\|_{L^2}^2+\frac12\|u^{n}_h\|_{L^2}^2)\notag.
	\end{align}
	
	For the second term on the right-hand side of \eqref{eq20180808_7}, using the Cauchy-Schwarz inequality, we get
	\begin{align}\label{eq20180808_8}
		&C(1+\|u_h^n\|_{L^2}^2)({\Delta} W_{n})^2(\|u^{n+1}_h\|_{L^2}^2+\frac12\|u^{n}_h\|_{L^2}^2)\\
		&=C(1+\|u_h^n\|_{L^2}^2)({\Delta} W_{n})^2(\|u^{n+1}_h\|_{L^2}^2-\| u^{n}_h\|_{L^2}^2+\frac32\|u^{n}_h\|_{L^2}^2)\notag\\
		&\le\theta_2\big(\|u^{n+1}_h\|_{L^2}^2-\|u^n_h\|_{L^2}^2)^2+(C+C\|u_h^n\|_{L^2}^4)({\Delta} W_{n})^4\notag\\
		&\quad+C\|u_h^n\|_{L^2}^4({\Delta} W_{n} \big)^2+C\|u_h^n\|_{L^2}^2({\Delta} W_{n})^2,\notag
	\end{align}
	where $\theta_2>0$ will be determined later.
	Using \eqref{lineargrow}, the third term on the right-hand side of \eqref{eq20180808_7} can be bounded by
	\begin{align}\label{eq20180808_9}
		&(G(u_h^n),u_h^n) {\Delta} W_{n}(\|u^{n+1}_h\|_{L^2}^2+\frac12\|u^{n}_h\|_{L^2}^2)\\
		&\quad=(G(u_h^n),u_h^n) {\Delta} W_{n}(\|u^{n+1}_h\|_{L^2}^2-\|u^{n}_h\|_{L^2}^2+\frac32\|u^{n}_h\|_{L^2}^2)\notag\\
		&\quad\le\theta_3(\|u^{n+1}_h\|_{L^2}^2-\|u^n_h\|_{L^2}^2)^2+(C+C\|u_h^n\|_{L^2}^4)({\Delta} W_{n})^2\notag\\
		&\qquad+\frac32(G(u_h^n),u_h^n)\|u_h^n\|_{L^2}^2{\Delta} W_{n}\notag,
	\end{align}
	where $\theta_3>0$ will be determined later.
	
	The third term on the right-hand side of \eqref{eq20180808_5} can be written as
	\begin{align}\label{eq20221227_4}
		&\frac12 \bigl(DG(u_h^n)\,G(u_h^n)((\Delta W_n)^2 - \tau),u_h^{n+1}\bigr)(\|u^{n+1}_h\|_{L^2}^2+\frac12\|u^{n}_h\|_{L^2}^2)\\
		&=\frac12 \bigl(DG(u_h^n)\,G(u_h^n)((\Delta W_n)^2 - \tau),u_h^{n+1}-u_h^n+u_h^n\bigr)(\|u^{n+1}_h\|_{L^2}^2+\frac12\|u^{n}_h\|_{L^2}^2)\notag\\
		&\le(\frac14\|u_h^{n+1}-u_h^n\|_{L^2}^2+C(1+\|u_h^n\|_{L^2}^2)((\Delta W_n)^2 - \tau)^2\notag\\
		&\quad+\frac12 \bigl(DG(u_h^n)\,G(u_h^n)((\Delta W_n)^2 - \tau),u_h^n\bigr)(\|u^{n+1}_h\|_{L^2}^2+\frac12\|u^{n}_h\|_{L^2}^2)\notag.
	\end{align}
	
	For the second term on the right-hand side of \eqref{eq20221227_4}, using the Cauchy-Schwarz inequality, we get
	\begin{align}\label{eq20221227_5}
		&C(1+\|u_h^n\|_{L^2}^2)((\Delta W_n)^2 - \tau)^2(\|u^{n+1}_h\|_{L^2}^2+\frac12\|u^{n}_h\|_{L^2}^2)\\
		&=C(1+\|u_h^n\|_{L^2}^2)((\Delta W_n)^2 - \tau)^2(\|u^{n+1}_h\|_{L^2}^2-\| u^{n}_h\|_{L^2}^2+\frac32\|u^{n}_h\|_{L^2}^2)\notag\\
		&\le\theta_2\big(\|u^{n+1}_h\|_{L^2}^2-\|u^n_h\|_{L^2}^2)^2+(C+C\|u_h^n\|_{L^2}^4)((\Delta W_n)^2 - \tau)^2\notag\\
		&\quad+C\|u_h^n\|_{L^2}^4((\Delta W_n)^2 - \tau)^2+C\|u_h^n\|_{L^2}^2((\Delta W_n)^2 - \tau)^2,\notag
	\end{align}
	where $\theta_4>0$ will be determined later.
	Using \eqref{lineargrow}, the third term on the right-hand side of \eqref{eq20221227_4} can be bounded by
	\begin{align}\label{eq20221227_6}
		&\frac12 \bigl(DG(u_h^n)\,G(u_h^n)((\Delta W_n)^2 - \tau),u_h^n)(\|u^{n+1}_h\|_{L^2}^2+\frac12\|u^{n}_h\|_{L^2}^2)\\
		&\quad=\frac12 \bigl(DG(u_h^n)\,G(u_h^n)((\Delta W_n)^2 - \tau),u_h^n)(\|u^{n+1}_h\|_{L^2}^2-\|u^{n}_h\|_{L^2}^2+\frac32\|u^{n}_h\|_{L^2}^2)\notag\\
		&\quad\le\theta_5(\|u^{n+1}_h\|_{L^2}^2-\|u^n_h\|_{L^2}^2)^2+(C+C\|u_h^n\|_{L^2}^4)(\Delta W_n)^2 - \tau)^2\notag\\
		&\qquad+\frac12 \bigl(DG(u_h^n)\,G(u_h^n)((\Delta W_n)^2 - \tau),u_h^n)\frac32\|u^{n}_h\|_{L^2}^2\notag,
	\end{align}
	where $\theta_5>0$ will be determined later.
	
	Choosing $\theta_1\sim\theta_5$ such that $\theta_1+\cdots+\theta_3\le\frac{1}{16}$, and then taking the summation over $n$ from $0$ to $\ell-1$ and taking the expectation on both sides of \eqref{eq20180808_5}, we obtain
	\begin{align}\label{eq20180808_10}
		&\frac38\E\left[\|u^{\ell}_h\|_{L^2}^4\right]+\frac{1}{16}\sum_{n=0}^{\ell-1}\E\left[(\|u^{n+1}_h\|_{L^2}^2-\|u^n_h\|_{L^2}^2)^2\right]+\sum_{n=0}^{\ell-1}\E\bigl[(\frac14\|(u^{n+1}_h-u^{n}_h)\|_{L^2}^2\\
		&\quad+\tau\|\nabla u_h^{n+1}\|_{L^2}^2+\frac{\tau}{2}\|u_h^{n+1}\|_{L^{q+1}}^{q+1})(\|u^{n+1}_h\|_{L^2}^2+\frac12\|u^{n}_h\|_{L^2}^2)\bigr]\notag\\
		&\le C\tau\sum_{n=0}^{\ell-1}\E\left[\|u^{n+1}_h\|_{L^2}^4\right]+C\tau\sum_{n=0}^{\ell-1}\E\left[\|\nabla u^{n+1}_h\|_{L^2}^{2(q+1)}\right]+\frac38\E\left[\|u^0_h\|_{L^2}^4\right]\notag\\
		&\quad+C\tau\sum_{n=0}^{\ell-1}\E\left[\|u_h^n\|_{L^2}^4\right]+C.\notag
	\end{align}
	
	When $\tau\le C$, we have
	\begin{align}\label{eq20180808_11}
		&\frac14\E\left[\|u^{\ell}_h\|_{L^2}^4\right]+\frac{1}{16}\sum_{n=0}^{\ell-1}\E\left[(\|u^{n+1}_h\|_{L^2}^2-\|u^n_h\|_{L^2}^2)^2\right]+\sum_{n=0}^{\ell-1}\E\bigl[(\frac14\|(u^{n+1}_h-u^{n}_h)\|_{L^2}^2\\
		&\quad+\tau\|\nabla u_h^{n+1}\|_{L^2}^2+\frac{\tau}{2}\|u_h^{n+1}\|_{L^4}^4)(\|u^{n+1}_h\|_{L^2}^2+\frac12\|u^{n}_h\|_{L^2}^2)\bigr]\notag\\
		&\le C\tau\sum_{n=0}^{\ell-1}\E\left[\|u^{n}_h\|_{L^2}^4\right]+C\tau\sum_{n=0}^{\ell-1}\E\left[\|\nabla u^{n+1}_h\|_{L^2}^{2(q+1)}\right]+\frac38\E\left[\|u^0_h\|_{L^2}^4\right]+C\notag.
	\end{align}
	
	Using Gronwall's inequality, we obtain
	\begin{align}\label{eq20180808_12}
		&\frac14\E\left[\|u^{\ell}_h\|_{L^2}^4\right]+\frac{1}{16}\sum_{n=0}^{\ell-1}\E\left[(\|u^{n+1}_h\|_{L^2}^2-\|u^n_h\|_{L^2}^2)^2\right]\\
		&\qquad+\sum_{n=0}^{\ell-1}\E\bigg[(\frac14\|(u^{n+1}_h-u^{n}_h)\|_{L^2}^2+\tau\|\nabla u_h^{n+1}\|_{L^2}^2\notag\\
		&\qquad+\frac{\tau}{2}\|u_h^{n+1}\|_{L^4}^4)(\|u^{n+1}_h\|_{L^2}^2+\frac12\|u^{n}_h\|_{L^2}^2)\bigg]\le C.\notag
	\end{align}
	
	\smallskip
	{\rm Step 2.} Similar to Step 1, using \eqref{eq20180808_5}--\eqref{eq20180808_9}, we have
	\begin{align}\label{eq20180808}
		&\frac38(\|u^{n+1}_h\|_{L^2}^4-\|u^n_h\|_{L^2}^4)+\frac{1}{16}(\|u^{n+1}_h\|_{L^2}^2-\|u^n_h\|_{L^2}^2)^2\\
		&+(\frac14\|(u^{n+1}_h-u^{n}_h)\|_{L^2}^2+\tau\|\nabla u_h^{n+1}\|_{L^2}^2+\frac{\tau}{2}\|u_h^{n+1}\|_{L^4}^4)(\|u^{n+1}_h\|_{L^2}^2+\frac12\|u^{n}_h\|_{L^2}^2)\notag\\
		&\le C\tau\|u^{n+1}_h\|_{L^2}^4+C\tau\|\nabla u_h^{n+1}\|_{L^2}^{2(q+1)}+(C+C\|u_h^n\|_{L^2}^4)({\Delta} W_{n})^4\notag\\
		&+(C+C\|u_h^n\|_{L^2}^4)({\Delta} W_{n})^2+(G(u_h^n),u_h^n)\|u_h^n\|_{L^2}^2{\Delta} W_{n}.\notag
	\end{align}
	
	Similar to Step 1, multiplying \eqref{eq20180808} by $\|u^{n+1}_h\|_{L^2}^4+\frac12\|u^{n}_h\|_{L^2}^4$,
	we can obtain the 8-th moment of the $L^2$ stability result of the discrete solution. Then repeating this
	process, the second moment of the $L^2$ stability result of the discrete solution can be obtained.
	
	\smallskip
	{\rm Step 3.} Suppose $2^{r-1}\le p\le 2^r$, and then by Young's inequality, we have
	\begin{align}
		\E\left[\|u^{\ell}_h\|_{L^2}^p\right]&\le \E\left[\|u^{\ell}_h\|_{L^2}^{2^r}\right]+C\le C,
	\end{align}
	where Step 2 is used in the second inequality. The proof is complete.
\end{proof}

\subsection{Error estimates of the finite element approximation}\label{subsec3}
In this subsection, we consider error estimates between the semi-discrete solution $u^n$ of Algorithm 1 and its finite element approximation $u^n_h$ from Algorithm 2. Let $e_h^n=u^n-u_h^n$ $(n = 0,1,2,\ldots,N)$.
In the following theorem, the $L^2$-projection is used in the proof of the error estimates and
the strong convergence rate is given.

\begin{theorem}\label{thm:derrest}
	Let $\{u^n\}$ and $\{ u_h^n \}_{n=1}^N$ denote respectively the solutions of Algorithm 1 and Algorithm 2. Then, under the condition \eqref{eq20180907}, there holds
	\begin{align*}
		&\sup_{0 \leq n \leq N} \E \left[\| e_h^n \|^2_{L^2} \right]
		+ \E \left[\tau\sum_{n=1}^N  \|\nabla e_h^n\|^2_{L^2}  \right]\le  Ch^2|\ln h|^{2(q-1)}.
	\end{align*}
\end{theorem}
\begin{proof}
	We write $e_h^n = \eta^n + \xi^n$ where
	\begin{align*}
		\eta^n := u^n - P_h u^n \quad \text{and} \quad \xi^n := P_h u^n - u_h^n,  \quad n = 0,1,2,...,N.
	\end{align*}

	Subtracting \eqref{dfem} from \eqref{milsteinscheme1} and setting $v_h = \xi^{n+1}$, the following error equation holds $\P$-almost surely,
	\begin{align} \label{eq20230205_1}
		&(\xi^{n+1} - \xi^n, \xi^{n+1}) = -(\eta^{n+1} - \eta^n, \xi^{n+1})
		- \tau ( \nab u^{n+1}-\nabla u^{n+1}_h, \nabla \xi^{n+1} ) \\
		&\qquad +  \tau\bigl(F(u^{n+1})-I_hF^{n+1}, \xi^{n+1}\bigr)+ (G(u^n)-G(u^n_h), \xi^{n+1}) \, \Delta W_n, \notag \\
		&\qquad +\frac12 \bigl((DG(u^n)\,G(u^n)-DG(u_h^n)\,G(u_h^n))((\Delta W_n)^2 - \tau),\xi^{n+1}\bigr)\notag\\
		& := T_1 + T_2 + T_3 + T_4 +T_5. \notag
	\end{align}
	
	The expectation of the left-hand side of \eqref{eq20230205_1} can be bounded by
	\begin{align} \label{derrest:4}
		\E \bigl[(\xi^{n+1} - \xi^n, \xi^{n+1}) \bigr]
		&= \frac{1}{2} \E \bigl[ \|\xi^{n+1}\|_{L^2}^2 - \|\xi^{n}\|_{L^2}^2 \bigr]  +  \frac{1}{2} \E \bigl[ \| \xi^{n+1} - \xi^n \|^2_{L^2} \bigr]. 
	\end{align}
	
	The first term on the right-hand side of \eqref{eq20230205_1} is 0 by the property of the $L^2$-projection.
	
	For the second term on the right-hand side of \eqref{eq20230205_1}, we have
	\begin{align} \label{eq20230205_2}
		\E \left[ T_2 \right] &= -\tau \E\left[(\nabla \eta^{n+1} + \nabla \xi^{n+1}, \nabla \xi^{n+1}\right] \\
		&\leq C\tau\E \left[\| \nabla \eta^{n+1} \|^2_{L^2} \right]- \frac{3}{4} \tau\E \left[ \| \nabla \xi^{n+1} \|^2_{L^2} \right]  \notag\\
		&\leq C\tau h^2 \mE[\|u^{n+1}\|^2_{H^2}]
		- \frac{3}{4} \E \left[\| \nabla \xi^{n+1} \|^2_{L^2} \right] \tau. \notag
	\end{align}
	
	In order to estimate the third term on the right-hand side of \eqref{eq20230205_1}, we write
	\begin{align}\label{eq20230205_3}
		\tau	\bigl( F(u^{n+1}) - I_hF^{n+1}, \xi^{n+1} \bigr) &= \tau\bigl( F(u^{n+1}) - F(P_h u^{n+1}), \xi^{n+1} \bigr) \\
		&\quad +  \tau\bigl( F(P_h u^{n+1}) - F^{n+1}, \xi^{n+1} \bigr) \notag\\
		&\quad +  \tau\bigl( F^{n+1}-I_hF^{n+1}, \xi^{n+1} \bigr) \notag.
	\end{align}
	
	{ The first term on the right-hand side of \eqref{eq20230205_3} can be bounded as follows. Using Cauchy-Schwarz's inequality, the Ladyzhenskaya inequality $\|u\|_{L^4} \leq C \|u\|^{1/2}_{L^2}\|\nab u\|^{1/2}_{L^2}$, and \eqref{Ph1} we obtain
		\begin{align}\label{eq4.57}
			&\tau\bigl( F(u^{n+1}) - F(P_h u^{n+1}), \xi^{n+1} \bigr)\\\nonumber
			&= -\tau\Bigl(\eta^{n+1}\Bigl[\sum_{i=0}^{q-1} (u^{n+1})^i (P_hu^{n+1})^{q-1-i} - 1\Bigr], \xi^{n+1}\Bigr) \\\nonumber
			&\leq \tau \|\eta^{n+1}\|_{L^4} \Bigl\|\sum_{i=0}^{q-1} (u^{n+1})^i (P_hu^{n+1})^{q-1-i} - 1\Bigr\|_{L^4} \|\xi^{n+1}\|_{L^2}\\\nonumber
			&\leq C\tau \|\eta^{n+1}\|^{1/2}_{L^2}\|\nab \eta^{n+1}\|^{1/2}_{L^2} \Bigl\|\sum_{i=0}^{q-1} (u^{n+1})^i (P_hu^{n+1})^{q-1-i} - 1\Bigr\|_{L^4} \|\xi^{n+1}\|_{L^2}\\\nonumber
			&\leq C\tau h \|\nab u^{n+1}\|^{1/2}_{L^2}\|\Delta u^{n+1}\|^{1/2}_{L^2} \Bigl\|\sum_{i=0}^{q-1} (u^{n+1})^i (P_hu^{n+1})^{q-1-i} - 1\Bigr\|_{L^4} \|\xi^{n+1}\|_{L^2}\\\nonumber
			&\leq C\tau h^2 \|\nab u^{n+1}\|_{L^2}\|\Delta u^{n+1}\|_{L^2} \Bigl\|\sum_{i=0}^{q-1} (u^{n+1})^i (P_hu^{n+1})^{q-1-i} - 1\Bigr\|^2_{L^4} \\\nonumber&\qquad\qquad\qquad+ \tau \|\xi^{n+1}\|^2_{L^2}.
		\end{align}
		
		Taking the summation $\sum_{n=0}^{\ell}$ to \eqref{eq4.57} for any $0\leq\ell \leq N-1$ we obtain
		\begin{align}\label{eq4.58}
			&\tau\sum_{n = 0}^{\ell}\bigl( F(u^{n+1}) - F(P_h u^{n+1}), \xi^{n+1} \bigr) \\\nonumber
			&\leq Ch^2\tau\sum_{n=0}^{\ell} \|\nab u^{n+1}\|_{L^2}\|\Delta u^{n+1}\|_{L^2} \Bigl\|\sum_{i=0}^{q-1} (u^{n+1})^i (P_hu^{n+1})^{q-1-i} - 1\Bigr\|^2_{L^4} \\\nonumber&\qquad\qquad\qquad+ \tau\sum_{n=0}^{\ell} \|\xi^{n+1}\|^2_{L^2}\\\nonumber
			&\leq Ch^2\Bigl(\tau\sum_{n=0}^{\ell} \|\nab u^{n+1}\|^2_{L^2}\|\Delta u^{n+1}\|^2_{L^2}\Bigr)^{1/2}\\\nonumber
			&\qquad\times\Bigl(\tau\sum_{n=0}^{\ell}\Bigl\|\sum_{i=0}^{q-1} (u^{n+1})^i (P_hu^{n+1})^{q-1-i} - 1\Bigr\|_{L^4}^4 \Bigr)^{1/2}  + \tau\sum_{n=0}^{\ell} \|\xi^{n+1}\|^2_{L^2}.
		\end{align}

		Next, applying the expectation to \eqref{eq4.58} and using Cauchy-Schwarz's inequality, and then using Lemma \ref{lemma3.3} we have
		\begin{align}
			&\mE\Bigl[\tau\sum_{n = 0}^{\ell}\bigl( F(u^{n+1}) - F(P_h u^{n+1}), \xi^{n+1} \bigr)\Bigr]\\\nonumber
			&\leq Ch^2\mE\Bigl[\Bigl(\tau\sum_{n=0}^{\ell} \|\nab u^{n+1}\|^2_{L^2}\|\Delta u^{n+1}\|^2_{L^2}\Bigr)^{1/2}\\\nonumber
			&\qquad\times\Bigl(\tau\sum_{n=0}^{\ell}\Bigl\|\sum_{i=0}^{q-1} (u^{n+1})^i (P_hu^{n+1})^{q-1-i} - 1\Bigr\|_{L^4}^4 \Bigr)^{1/2} \Bigr] +\mE\Bigl[ \tau\sum_{n=0}^{\ell} \|\xi^{n+1}\|^2_{L^2}\Bigr]\\\nonumber
			&\leq Ch^2\Bigl(\mE\Bigl[\tau\sum_{n=0}^{\ell}\|\nab u^{n+1}\|^2_{L^2}\|\Delta u^{n+1}\|^2_{L^2}\Bigr]\Bigr)^{1/2}\\\nonumber
			&\qquad\times \Bigl(\mE\Bigl[\tau\sum_{n=0}^{\ell}\Bigl\|\sum_{i=0}^{q-1} (u^{n+1})^i (P_hu^{n+1})^{q-1-i} - 1\Bigr\|_{L^4}^4\Bigr]\Bigr)^{1/2}\\\nonumber&\qquad\qquad\qquad+\mE\Bigl[ \tau\sum_{n=0}^{\ell} \|\xi^{n+1}\|^2_{L^2}\Bigr]\\\nonumber
			&\leq Ch^2 \Bigl(\mE\Bigl[\tau\sum_{n=0}^{\ell}\Bigl\|\sum_{i=0}^{q-1} (u^{n+1})^i (P_hu^{n+1})^{q-1-i} - 1\Bigr\|_{L^4}^4\Bigr]\Bigr)^{1/2}\\\nonumber&\qquad\qquad\qquad+\mE\Bigl[ \tau\sum_{n=0}^{\ell} \|\xi^{n+1}\|^2_{L^2}\Bigr].
		\end{align}
		
		Moreover, using the embedding inequality $\|u\|_{L^r} \leq C\|u\|_{H^1}$ for any integers $r \geq 2$ (see \cite[Corollary 9.14]{Brezis}) we also have
		\begin{align*}
			&\mE\Bigl[\tau\sum_{n=0}^{\ell}\Bigl\|\sum_{i=0}^{q-1} (u^{n+1})^i (P_hu^{n+1})^{q-1-i} - 1\Bigr\|_{L^4}^4\Bigr] \\\nonumber
			&\leq C\mE\Bigl[\tau\sum_{n=0}^{\ell}\Bigl(\|u^{n+1}\|^{4(q-1)}_{L^{4(q-1)}} + \|P_h u^{n+1}\|^{4(q-1)}_{L^{4(q-1)}}+ C\Bigr)\Bigr]\\\nonumber
			&\leq C\mE\Bigl[\tau\sum_{n=0}^{\ell}\Bigl(\| u^{n+1}\|^{4(q-1)}_{H^1} + \|P_h u^{n+1}\|^{4(q-1)}_{H^1}+ C\Bigr)\Bigr]\leq C,
		\end{align*}
		where the last inequality is obtained by using Lemma \ref{lemma3.3}. In summary, we obtain the following estimate for the first term of $T_3$
		\begin{align}\label{T_31}
			\mE\Bigl[\tau\sum_{n = 0}^{\ell}\bigl( F(u^{n+1}) - F(P_h u^{n+1}), \xi^{n+1} \bigr)\Bigr] \leq Ch^2 + \mE\Bigl[ \tau\sum_{n=0}^{\ell} \|\xi^{n+1}\|^2_{L^2}\Bigr].
		\end{align}

	}
	
	By using the one-sided Lipchitz condition \eqref{oneside_Lip}, the second term on the right-hand side of \eqref{eq20230205_3} can be bounded by
	\begin{align}
		\label{derrest:9_1}
		\E \left[\bigl( F(P_h u^{n+1}) - F^{n+1}, \xi^{n+1} \bigr) \right] \leq \E \left[\| \xi^{n+1} \|^2_{L^2} \right].
	\end{align}
	
	Using properties of the interpolation operator, the inverse inequality, and the fact that $u_h^{n+1}$ is a piecewise linear polynomial, the third term on the right-hand side of \eqref{eq20230205_3} can be handled by
	\begin{align}\label{eq20180711_17}
		&\E \left[\bigl( F^{n+1}-I_hF^{n+1}, \xi^{n+1} \bigr) \right]\\
		&\quad\le \E \Bigl[Ch^2\sum_{K\in\mathcal{T}_h}\|q(u_h^{n+1})^{q-1}\nabla u_h^{n+1}\|^2_{L^2(K)} \Bigr]+\E \left[\| \xi^{n+1} \|^2_{L^2} \right]\notag\\
		&\quad\le \E \Bigl[Ch^2\left(\| u_h^{n+1}\|_{L^\infty}^{2(q-1)}\|\nabla u_h^{n+1} \|^2_{L^2}\right) \Bigr]+\E \left[\| \xi^{n+1} \|^2_{L^2} \right]\notag\\
		&\quad\le \E \Bigl[Ch^2|\ln h|^{2(q-1)}\Bigl(\sum_{K\in\mathcal{T}_h}(\|\nabla u_h^{n+1} \|^2_{L^2(K)}+\|u_h^{n+1} \|^2_{L^2(K)})\Bigr)^{q-1}\notag\\
		&\qquad \quad\|\nabla u_h^{n+1} \|^2_{L^2} \Bigr]+\E \left[\| \xi^{n+1} \|^2_{L^2} \right]\notag
		\\
		&\quad\le \E \left[Ch^2|\ln h|^{2(q-1)}(\|u_h^{n+1} \|^{2(q-1)}_{L^2}+\|\nabla u_h^{n+1} \|^{2(q-1)}_{L^2})\|\nabla u_h^{n+1} \|^2_{L^2} \right]\notag\\
		&\qquad\quad+\E \left[\| \xi^{n+1} \|^2_{L^2} \right]\notag\\
		&\quad\le \E \left[Ch^2|\ln h|^{2(q-1)}(\|u_h^{n+1} \|^{2q}_{L^2}+\|\nabla u_h^{n+1} \|^{2q}_{L^2}) \right]+\E \left[\| \xi^{n+1} \|^2_{L^2} \right]\notag\\
		&\quad\le Ch^2|\ln h|^{2(q-1)}+\E \left[\| \xi^{n+1} \|^2_{L^2} \right]\notag.
	\end{align}
	
	{
		Combining \eqref{T_31}--\eqref{eq20180711_17} yields
		\begin{align} \label{derrest:11}
			\E \left[ \tau \sum_{n=0}^{\ell} T_3 \right] \leq Ch^2|\ln h|^{2(q-1)}+C \E \left[\tau\sum_{n=0}^{\ell}\|\xi^{n+1}\|_{L^2}^2 \right] .
		\end{align}
	}
	
	By the assumption ${\bf (A1)}$ for $G(\cdot)$ and then Lemma \ref{lemma3.3}, we have
	\begin{align}\label{derrest:15}
		\E [T_4] &\leq \frac{1}{2} \E \left[\|\xi^{n+1} - \xi^n\|^2_{L^2} \right]+ \frac{1}{2}\tau \E \Bigl[\|G(u^n)-G(u^n_h)\|^2_{L^2} \, ds \Bigr]  \\
		&\leq \frac{1}{2} \E \left[\|\xi^{n+1} - \xi^n\|^2_{L^2} \right]+ \frac{1}{2}\tau \E \Bigl[\|u^n-u^n_h\|^2_{L^2} \, ds \Bigr] \notag \\
		&\leq \frac{1}{2} \E \left[\|\xi^{n+1} - \xi^n\|^2_{L^2} \right] + C \E \left[\| \eta^n +  \xi^n\|^2_{L^2} \right] \tau \notag \\
		&\leq \frac{1}{2} \E \left[\|\xi^{n+1} - \xi^n\|^2_{L^2} \right] + C\E \left[\|  \xi^{n} \|^2_{L^2} \right] \tau + C\tau h^2 \mE\bigl[\|\nab u^n\|^2_{L^2}\bigr]\notag\\
		&\leq \frac{1}{2} \E \left[\|\xi^{n+1} - \xi^n\|^2_{L^2} \right] + C\E \left[\|  \xi^{n} \|^2_{L^2} \right] \tau + C\tau h^2.\notag
	\end{align}
	
	By using the assumption {\bf (A3)} for $G$ and then Lemma \ref{lemma3.3}, we have
	\begin{align}\label{eq20220205_3}
		\E [T_5] &\leq \frac{1}{2} \E \left[\|\xi^{n+1} - \xi^n\|^2_{L^2} \right]+ \frac{1}{2}\E \Bigl[\|DG(u^n)\,G(u^n)-DG(u^n)\\
		&\quad\cdot G(u_h^n)\|^2_{L^2}\Bigr]\tau^2  + \frac{1}{2}\E \Bigl[\|DG(u^n)\,G(u_h^n)-DG(u_h^n)\,G(u_h^n)\|^2_{L^2}\Bigr]\tau^2\notag\\
		&\leq \frac{1}{2} \E \left[\|\xi^{n+1} - \xi^n\|^2_{L^2} \right]+ C \E \Bigl[\|u^n-u^n_h\|^2_{L^2} \, ds \Bigr]\tau^2 \notag \\
		&\leq \frac{1}{2} \E \left[\|\xi^{n+1} - \xi^n\|^2_{L^2} \right] + C \E \left[\| \eta^n +  \xi^n\|^2_{L^2} \right] \tau^2 \notag \\
		&\leq \frac{1}{2} \E \left[\|\xi^{n+1} - \xi^n\|^2_{L^2} \right] + C\E \left[\|  \xi^{n} \|^2_{L^2} \right] \tau^2 + C\tau^2 h^2\mE\bigl[\|\nab u^n\|^2_{L^2}\bigr]\notag\\
		&\leq \frac{1}{2} \E \left[\|\xi^{n+1} - \xi^n\|^2_{L^2} \right] + C\E \left[\|  \xi^{n} \|^2_{L^2} \right] \tau^2 + C\tau^2 h^2\notag
	\end{align}
	
	Taking the expectation on \eqref{eq20230205_1} and combining estimates \eqref{derrest:4}--\eqref{eq20220205_3}, summing over
	$n = 0, 1, 2, ..., \ell-1$ with $1 \leq \ell \leq N$, and using Lemma \ref{lemma3.3} we obtain
	\begin{align} \label{derrest:17}
		\frac14&\E \left[\| \xi^{\ell} \|^2_{L^2} \right]
		+ \frac{1}{4} \E \Bigl[\tau \sum_{n = 1}^{\ell} \| \nabla \xi^n \|^2_{L^2} \Bigr] \nonumber\\
		& \leq \frac{1}{2} \E \left[ \| \xi^0 \|^2_{L^2} \right] + C\E \Bigl[ \tau \sum_{n=0}^{\ell-1} \| \xi^n \|^2_{L^2} \Bigr] \\\nonumber
		&\qquad+Ch^2|\ln h|^{2(q-1)} + Ch^2\mE\Bigl[\tau\sum_{n=0}^{\ell-1}\|u^{n+1}\|^2_{H^2}\Bigr]\notag\\
		& \leq \frac{1}{2} \E \left[ \| \xi^0 \|^2_{L^2} \right] + C\E \Bigl[ \tau \sum_{n=0}^{\ell-1} \| \xi^n \|^2_{L^2} \Bigr] +Ch^2|\ln h|^{2(q-1)}.\notag
	\end{align}
	
	Finally, the conclusion of the theorem follows from the discrete Gronwall's inequality, the fact that $\xi^0 = 0$, and the triangle inequality.
\end{proof}

\subsection{Global error estimates}
Finally, we are ready to state the global error estimates of our proposed method in the following theorem. 
\begin{theorem}\label{thm:global}
	Let $u$ and $\{ u_h^n \}_{n=1}^N$ denote respectively the solutions of \eqref{weak_form} and Algorithm 2. Then, under the conditions of Theorem \ref{theorem_semi} and Theorem \ref{thm:derrest}, there holds
	\begin{align*}
		\sup_{0 \leq n \leq N} &\E \left[\| u(t_n) - u^n_h\|^2_{L^2} \right]
		\\\nonumber
		&\qquad+ \E \left[\tau\sum_{n=1}^N  \|\nabla(u(t_n) - u^n_h)\|^2_{L^2}  \right] \le  C\bigl(\tau^{2(1-\eps)} +  h^2|\ln h|^{2(q-1)}\bigr).
	\end{align*}
\end{theorem}

\section{Numerical Experiments}\label{nume}
In this section, three numerical tests are presented. In Test 1, the evolution and stability of \eqref{eq1.1} in the case $F(u) = u-u^3$ are illustrated with different noise intensities. Test 2 provides the visualization of the stability using a different drift term and diffusion term. Test 3 presents the error orders with respect to time step size $\tau$. {The domain $\D$ for all the following tests is chosen to be $\D=[-1,1]\times[-1,1]$.}\\

{\bf Test 1.} Consider the initial condition:
\begin{equation}\label{IC}
	u_0(x,y)=\tanh(\frac{\sqrt{x^2+y^2}-0.6}{\sqrt{2}\epsilon}).
\end{equation}

\begin{figure}[h]
	\subfloat[$\delta=0.1$]{\includegraphics[scale=0.45]{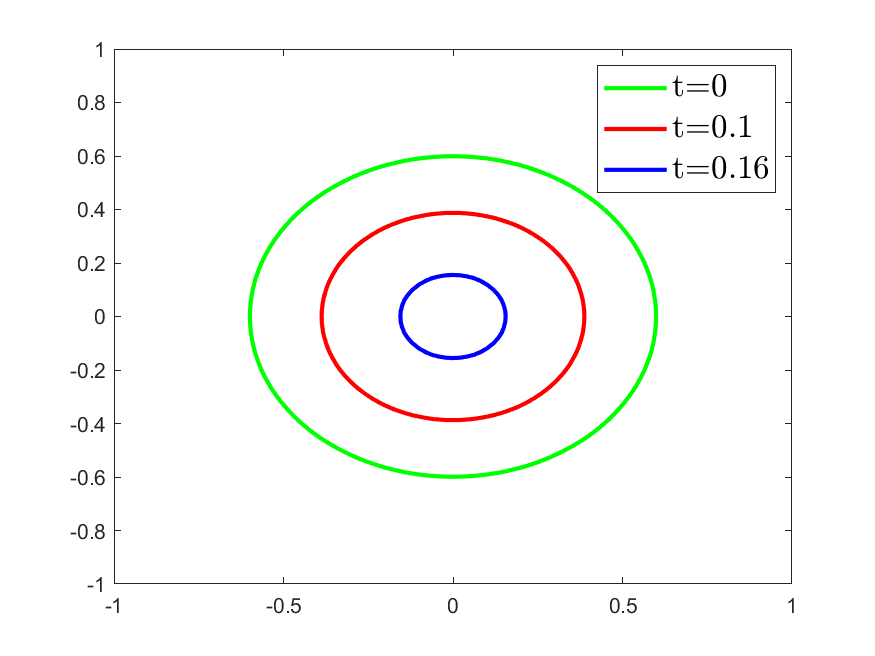}}
	\subfloat[$\delta=1$]{\includegraphics[scale=0.45]{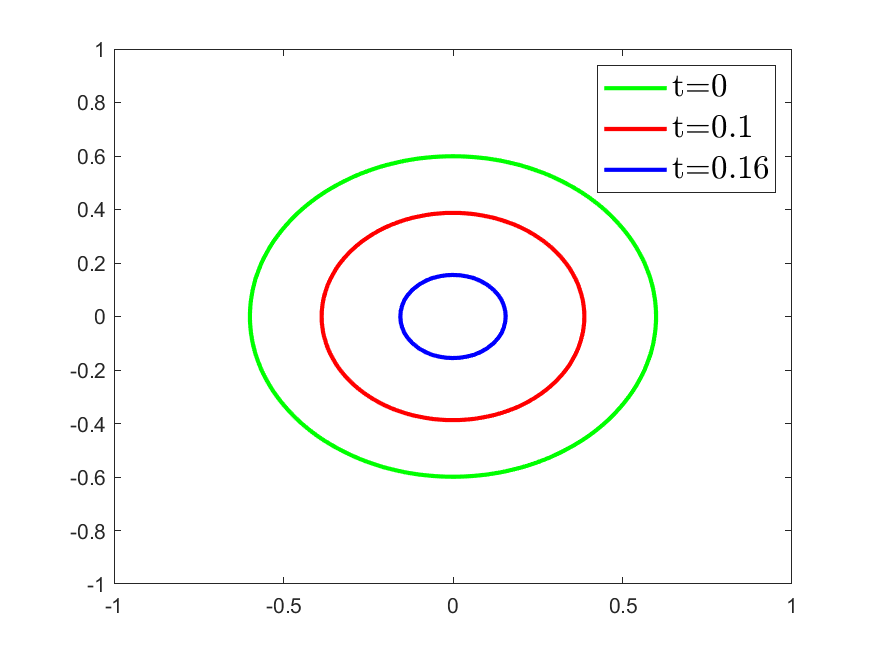}}
	\caption{{Test 1:} Zero-level sets of the evolution: $\tau=5\times 10^{-4}$, $h=0.02$, $\eps=0.04$.}      
	\label{test1-evo}
\end{figure}

\begin{figure}[h]
	\subfloat[$\delta=0.1$]{\includegraphics[scale=0.45]{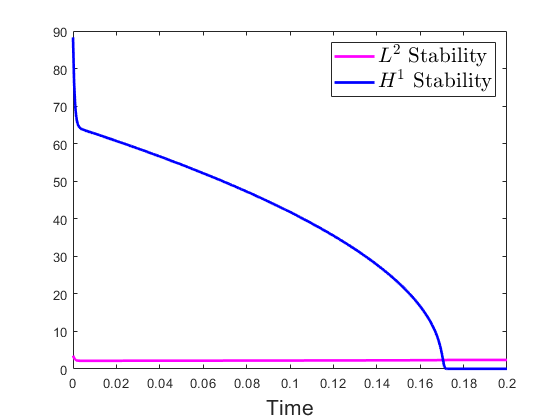}}
	\subfloat[$\delta=1$]{\includegraphics[scale=0.45]{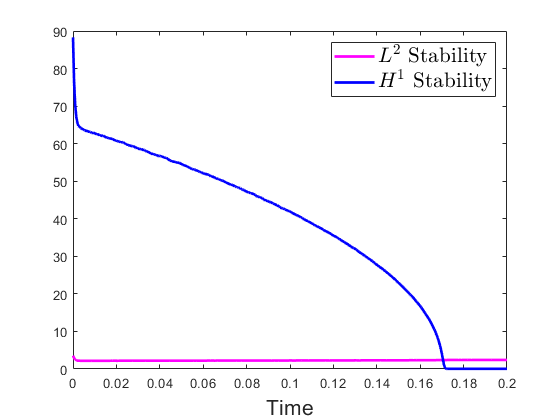}}
	\caption{{\color{blue}Test 1:} Stability demonstration (average): $\tau=5\times 10^{-4}$, $h=0.02$, $\eps=0.04$.}
	\label{test1-ave-stabi}
\end{figure}

\begin{figure}[h]
	\subfloat[$\delta=0.1$]{\includegraphics[scale=0.45]{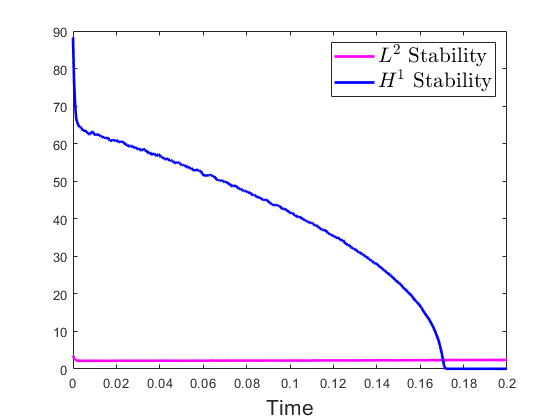}}
	\subfloat[$\delta=1$]{\includegraphics[scale=0.45]{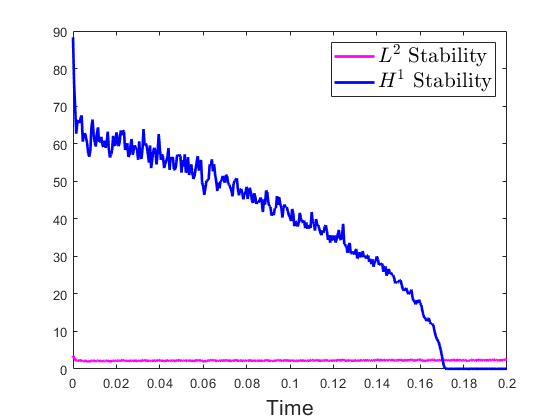}}
	\caption{{\color{blue}Test 1:} Stability demonstration (one sample point): $\tau=5\times 10^{-4}$, $h=0.02$, $\eps=0.04$.}
	\label{test1-one-stabi}
\end{figure}

For this test, $F(u)=u-u^3$ is used as the nonlinear term, and $G(u)=\delta u$ is used as the diffusion term. In Figure \ref{test1-evo}, the zero-level sets of the evolution using two different levels of noise intensity are shown. One can observe that the average zero-level set is a shrinking circle for both levels of noise intensity. Figure \ref{test1-ave-stabi} demonstrates the $\mathbb{E} L^2$ and $\mathbb{E} H^1$ stability for each time step. One can make the observation that they are both bounded. A one-sample $\mathbb{E} L^2$ and $\mathbb{E} H^1$ stability are provided in Figure \ref{test1-one-stabi}. Those stability results are still bounded but they are not always decreasing over time.   

\begin{figure}[h]
	\subfloat[$\delta=0.1$]{\includegraphics[scale=0.45]{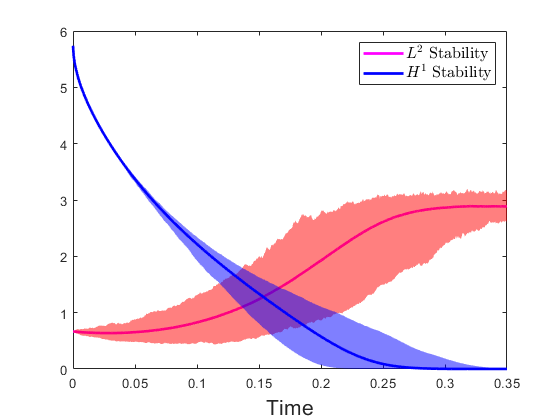}}
	\subfloat[$\delta=1$]{\includegraphics[scale=0.45]{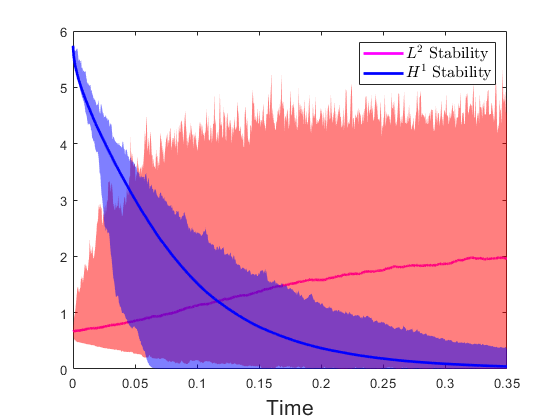}}
	\caption{{Test 2:} Stability demonstration (average and max/min): $\tau=5\times 10^{-4}$, $h=0.02$, $\eps=0.5$.}
	\label{test2-ave-stabi}
\end{figure}

{\bf Test 2}. For this test, the initial condition is still in (\ref{IC}), and that $\epsilon = 0.5$. The drift term is changed to $F(u)=u-u^{11}$, and the diffusion term is changed to $G(u)=\delta\sqrt{u^2+1}$. In Figure \ref{test2-ave-stabi}, the $\mathbb{E} L^2$ and $\mathbb{E} H^1$ stability are given by the blue and pink solid lines, along with the maximum and minimum of those two stabilities given by upper and lower edges of the shaded red and blue regions. One can see that both the $\mathbb{E} L^2$ and $\mathbb{E} H^1$ stability are bounded.

{\bf Test 3.}\label{Test 3} Consider the initial condition:
\begin{equation}\label{IC_2}
	u_0(x,y)=\tanh(\frac{\sqrt{x^2+y^2}-0.8}{\sqrt{2}\epsilon}).
\end{equation}


In this test, we use $\epsilon=0.3$, $F(u)=u-u^3$ as the drift term, and $G(u)=\delta u$ as the diffusion term. {The final time is $T=0.25$}. Table \ref{order_table} demonstrates the error $\{\underset{0\leq n \leq N}{\sup} \E[||e^n||^2_{L^2(\D)}]\}^{\frac{1}{2}}$ and the error $\{\E[\sum_{n=1}^N \tau ||\nabla e^n||^2_{L^2(\D)}]\}^{\frac{1}{2}}$. The error $\{\underset{0\leq n \leq N}{\sup} \E[||e^n||^2_{L^2(\D)}]\}^{\frac{1}{2}}$ is denoted by $L^\infty \E L^2$, and the error $\{\E[\sum_{n=1}^N \tau ||\nabla e^n||^2_{L^2(\D)}]\}^{\frac{1}{2}}$ is denoted by $\E L^2H^1$. By observingTable \ref{order_table}, one can see that the error orders for both $L^\infty \E L^2$ and $\E L^2H^1$ are 1.

\begin{table}[h]
	\begin{center}
		\begin{tabular}{ |c|c|c|c|c| } 
			\hline
			& $L^\infty \E L^2$ error & order & $\E L^2H^1$ error & order \\ 
			\hline
			$\tau=0.025$ & 0.080163 & - & 0.054115 & - \\ 
			\hline
			$\tau=0.0125$ & 0.038604 & 1.0542 & 0.027675 & 0.9675 \\ 
			\hline
			$\tau= 0.0625$ & 0.018036 & 1.0978 & 0.013978 & 0.9855 \\
			\hline
			$\tau= 0.03125$ & 0.008724 & 1.0479 & 0.007467 & 0.9045 \\
			\hline
		\end{tabular}
		\caption{\label{order_table} {Test 3:} Time step errors and rates of convergence of Test 3: $h = 0.01, \epsilon = 0.3, \delta = 0.01, {T=0.25}$.}
	\end{center}
\end{table}



\end{document}